\documentclass[11pt]{amsart}

\renewcommand{\bar}{\overline}

\newcommand{\pa}{\partial}
\newcommand{\ph}{\varphi}

\newcommand{\rr}{{|z| \leq\frac{\log m}{\sqrt m}}}  

\font\strange=msbm10

\newcommand{\C}{{{\mathchoice  {\hbox{$\textstyle{\text{\strange C}}$}}
{\hbox{$\textstyle{\text{\strange C}}$}}
{\hbox{$\scriptstyle{\text{\strange C}}$}}
{\hbox{$\scriptscriptstyle{\text{\strange C}}$}}}}}

\newcommand{\Z}{{{\mathchoice  {\hbox{$\textstyle{\text{\strange Z}}$}}
{\hbox{$\textstyle{\text{\strange Z}}$}}
{\hbox{$\scriptstyle  Z\kern-0.3em  Z$}}
{\hbox{$\scriptscriptstyle  Z\kern-0.2em  Z$}}}}}

\newcommand{\N}{{{\mathchoice  {\hbox{$\textstyle{\text{\strange N}}$}}
{\hbox{$\textstyle{\text{\strange N}}$}}
{\hbox{$\scriptstyle  N\kern-0.3em  N$}}
{\hbox{$\scriptscriptstyle  N\kern-0.2em  N$}}}}}

\newcommand{\oo}[1]{{O(\frac{1}{m^{{#1}}})}}
\newcommand{\ooo}[1]{{O(|z|^{{#1}})}}
\newcommand{\bb}{{\frac{\sqrt{-1}}{2\pi}}}

\newcommand{\crr}[4]{{R_{{#1}\bar{#2}{#3}\bar{#4}}}}
\newcommand{\css}[2]{{R_{{#1}\bar{#2}}}}

\newcommand{\gtt}[2]{{g^{{#1}\bar{#2}}}}

\newcommand{\gss}[3]{{\frac{\pa g_{{#1}\bar{#2}}}
{\pa z_{#3}}}}

\usepackage{amsmath,amsthm,amscd}

\newfont{\fnt}{cmr10 scaled 550}

\title
[
On the Lower Terms of the Asymptotic
Expansion]{On the Lower Order Terms of  the
Asymptotic Expansion of
Zelditch}
\author{Zhiqin Lu}
\date{Sept. 28, 1998}
\subjclass{Primary: 53A30; Secondary: 32C16}

\keywords{Szeg\"o kernel, asymptotic expansion, ample line bundle}
\address[Zhiqin Lu]
{Department of Mathematics\\
Columbia University\\
New York, NY 10027}
\email[Zhiqin Lu]{lu@cpw.math.columbia.edu}

\newtheorem{theorem}{Theorem}[section]
\newtheorem*{thm}{Theorem}
\newtheorem{lemma}{Lemma}[section]
\newtheorem{cor}{Corollary}[section]
\newtheorem{prop}{Proposition}[section]
\newtheorem{claim}{Claim}
\newtheorem{definition}{Definition}[section]

\newtheorem{example}{Example}
\theoremstyle{remark}
\newtheorem{rem}{Remark}[section]

\begin{document}
\maketitle

\numberwithin{equation}{section}

\tableofcontents

\section{Introduction}
A projective algebraic manifold $M$ is a complex manifold in certain
projective space $CP^m$, $m\geq\dim_{\C}M=n$. The hyperplane line bundle of
$CP^m$ restricts to an ample line bundle $L$ on $M$, which is called a
polarization on $M$. 
A K\"ahler metric $g$ is called a polarized metric,  if the corresponding
K\"ahler form  represents the first Chern class $c_1(L)$ of $L$ in
$H^2(M,\Z)$. Given any polarized K\"ahler metric $g$, there is a Hermitian
metric $h$ on $L$ whose Ricci form is equal to $\omega_g$. For
each positive integer $m>0$, the Hermitian metric $h$ induces a Hermitian
metric $h_m$ on $L^m$. Let  
$\{S_0^m,\cdots,S_{d_m-1}^m\}$ be an orthonormal basis of the space
$H^0(M,L^m)$ of all holomorphic
global sections of $L^m$. Such a basis $\{S_0^m,\cdots,S_{d_m-1}^m\}$
induces a holomorphic embedding $\ph_m$ of $M$ into $CP^{d_m-1}$ by
assigning the point $x$ of $M$ to $[S_0^m(x),\cdots,S_{d_m-1}^m(x)]$
in $CP^{d_m-1}$. Let $g_{{\fnt FS}}$ be the standard Fubini-Study metric
on
$CP^{d_m-1}$, i.e.,
$\omega_{{\fnt FS}}=\bb\pa\bar\pa\log\sum_{i=0}^{d_m-1}|w_i|^2$
for a homogeneous coordinate system $[w_0,\cdots,w_{d_m-1}]$ of
$CP^{d_m-1}$. The $\frac 1m$-multiple of $g_{{\fnt FS}}$ on $CP^{d_m-1}$
restricts
to a K\"ahler metric $\frac 1m\ph^*_mg_{{\fnt FS}}$ on $M$. This metric is
a
polarized K\"ahler metric on $M$ and is called the Bergman metric with
respect to $L$.

One of the main theorem in~\cite{T5} is the following

\begin{thm}[Tian] With all the notations as above, we have
\[
||\frac 1m\ph^*_mg_{{\fnt FS}}-g||_{C^2}=O(\frac{1}{\sqrt{m}})
\]
for any polarized metric on $M$ with respect to $L$.
\end{thm}

Using Tian's peak section method, Ruan~\cite{ru} proved the
$C^\infty$ convergence and
improved the bound to $\oo{}$.
Recently, S. Zelditch~\cite{sz} beautifully generalized 
the
above theorem by
using the Szeg\"o kernel on the unit circle bundle of $L^*$ over $M$. His
result gives the asymptotic expansion of the
potential of the  Bergman metric:

\begin{thm}[Zelditch]
Let $M$ be a compact complex manifold of dimension $n$ (over $\C$) and let
$(L,h)\rightarrow M$ be a positive Hermitian holomorphic line bundle. 
Let $x$ be a point of $M$.
Let
$g$ be the K\"ahler metric on $M$ corresponding to the K\"ahler form
$\omega_g=Ric(h)$. For each $m\in \N$, $h$ induces a Hermitian metric
$h_m$
on $L^m$. Let
$\{S_0^m,\cdots,S_{d_m-1}^m\}$ be any orthonormal basis of $H^0(M,L^m)$,
$d_m=\dim H^0(M,L^m)$, with respect to the inner product 
\[
(S_1,S_2)_{h_m}
=\int_Mh_m(S_1(x),S_2(x))dV_g,
\]
where $dV_g=\frac{1}{n!}\omega_g^n$ is the
volume form of $g$. Then there is a complete asymptotic expansion:
\begin{equation}\label{fud1}
\sum_{i=0}^{d_m-1}||S_i^m(x)||_{h_m}^2
=a_0(x)m^n+a_1(x)m^{n-1}+a_2(x)m^{n-2}+\cdots
\end{equation}
for certain smooth coefficients $a_j(x)$ with $a_0=1$. 
More precisely,
for
any $k$
\[
||\sum_{i=0}^{d_m-1}||S_i^m(x)||_{h_m}^2
-\sum_{j<R}a_j(x)m^{n-j}||_{C^k}\leq C_{R,k}m^{n-R},
\]
where $C_{R,k}$ depends on $R,k$ and the manifold $M$.
\end{thm}

\begin{rem}
The referee informed the author that Professor D. Catlin has also obtained
the above Zelditch's result independently.
\end{rem}

In this paper, we give a method to compute the coefficients $a_1(x),
a_2(x),\cdots$ ($a_0(x)=1$
was pointed out in~\cite{sz} in a more general setting). 
Our result is

\begin{theorem}\label{fud}
With the notations as in the above theorem,
each coefficient $a_j(x)$ is a polynomial of the curvature 
 and its covariant derivatives  at $x$ with  weight $j$
\footnote{See Definition~\ref{md1}.}. Such a
polynomial can be
found by finite many steps of algebraic operations. 
In particular,
\[
\left\{
\begin{array}{l}
a_0=1\\
a_1=\frac 12 \rho\\
a_2=\frac 13\Delta\rho+\frac{1}{24}(|R|^2-4|Ric|^2+3\rho^2)\\
a_3=\frac 18\Delta\Delta\rho+\frac 1{24}\,div\,div\, (R,Ric)
-\frac 16 div\,div (\rho Ric)\\
+\frac{1}{48}\Delta (|R|^2-4|Ric|^2+8\rho^2)
+\frac{1}{48}\rho(\rho^2-4|Ric|^2+|R|^2)\\
+\frac{1}{24}(\sigma_3(Ric)-Ric(R,R)-R(Ric,Ric)),
\end{array}
\right.
\]
where $R, Ric$ and $\rho$ represent the curvature tensor, the Ricci
curvature and the scalar curvature of $g$, respectively
and $\Delta$ represents the Laplace operator of $M$. For the precise
definition of the terms in the expression of $a_3$, see Section 5. 
\end{theorem}

We 
 use the peak section method initiated by Tian in~\cite{T5} to compute the
coefficients $a_j (j\in\N)$.
Consider $H^0(M,L^m)$ for $m$ large enough. Fixing a point $x\in M$, 
 by the standard $\bar\pa$-estimate Tian observed that the sections which
do not
vanish at $x$ at a very high order are known in the sense that
one can completely control their behavior around the point $x$. These
sections are called peak sections (in the terminology of ~\cite{T5}). We
proved that the coefficients $a_1,a_2,\cdots$ only depend on the inner
products of the  peak sections. Various techniques
are used to give
the
asymptotic expansion of these inner products, including some combinatorial
lemmas, to simplify the computation and thus make
the
computation feasible.

Zelditch's work is based on the analysis of the asymptotic expansion of
the Szeg\"o kernel on the unit circle bundle of the ample line bundle
over a complex manifold. To be more precise, let $C$ be the unit circle
bundle and let $\Pi(x,y)$ be its Szeg\"o kernel (with the natural
measure). Since $C$ is $S^1$ invariant, we have Fourier coefficients
\[
\Pi_m(x,x)=\int_{S^1} e^{-im\theta}\Pi(r_\theta x,x) d\theta,
\]
where $r_\theta$ is the circle action. The key observation
by Zelditch  is that
\[
\sum_{i=0}^{d_m-1}||S_i^m(x)||^2_{h_m}=\Pi_m(x,x).
\]
Thus the general theory of Szeg\"o kernel can be applied.

There are a lot of works on the Bergman and Szeg\"o kernels on the 
pseudoconvex domain on $\C^n$ (~\cite{BM1}~\cite{BS1},~\cite{BEG},
~\cite{G1},~\cite{KH1} and~\cite{KH2}, for example) following the program
of
Fefferman~\cite{F1}.
While our method is completely complex-geometric, it should also be 
possible to compute the coefficients from the general 
theory of Szeg\"o kernel. In particular, we noticed the works
in~\cite{BEG} and~\cite{KH1}, the coefficients are proved to be the Weyl
functionals of the curvature tensor of the ambient metric defined by
Fefferman. But I don't know how to relate it to the curvature of the base
manifold.

The organization of the paper is as follows: in \S 2 we introduce
the concept of peak sections initiated by Tian and discuss their
properties; in \S 3 we give the iteration process; in \S 4 we
prove the main theorem of this paper except for the computation of the
$a_3$ term; and in \S 5 and in the Appendix we calculate the
$a_3$
term.

{\bf Acknowledgment.}  
 The idea that peak section
method can be used to get the  result is from G. Tian.
The author thanks him for the encouragement and help during the
preparation of this paper. He also thanks 
S. Zelditch for sending him the preprint~\cite{sz} 
and bringing him to the attention of the works of K. Hirachi.
The author thanks  K. Hirachi,  L. Boutet de
Monvel and D. Phong 
for their help.
Finally, the author deeply thanks the referee for the careful 
proofreading and  many
suggestions
to improve the organization and the style of this paper.

\section{Peak Global Sections}
Let $M$ be an $n$-dimensional algebraic manifold with a
positive Hermitian line bundle $(L,h)\rightarrow M$. Suppose 
that the K\"ahler
form $\omega_g$ is defined by the curvature $Ric(h)$ of $h$. That is, 
fixing a point $x_0\in M$, under
local coordinate $(z_1,\cdots,z_n)$ at $x_0$,
\[
\omega_g=-\bb\sum_{\alpha,\beta=1}^n
\frac{\pa^2}{
\pa z_\alpha\pa\bar z_\beta}\log a\,
dz_\alpha\wedge d\bar z_\beta=\bb\sum_{\alpha,\beta=1}^n
g_{\alpha\bar\beta} dz_\alpha\wedge d\bar z_\beta,
\]
where $a$ is the local representation of  the Hermitian metric
$h$.

Let $d=d_m=\dim_{\C} H^0(M,L^m)$ for a fixed 
integer $m$. Let $S_0,\cdots, S_{d-1}$
be
a basis of $H^0(M,L^m)$. The metrics $(h,\omega_g)$ define an
inner product $(\,,\,)$ on $H^0(M,L^m)$ as 
\[
(S_A,S_B)_{h_m}=\int_M<S_A,S_B>_{h_m}dV_g,
\quad A,B=0,\cdots,d-1,
\]
where $<S_A,S_B>_{h_m}$ is the pointwise inner product with respect
to $h_m$ and $dV_g=\frac{1}{n!}\omega_g^n$. 

We assume that at the point $x_0\in M$,
\[
S_0(x_0)\neq 0,\quad S_{A}(x_0)=0,\quad A=1,\cdots, d-1.
\]
Suppose
\[
F_{AB}=(S_A,S_B)_{h_m}, \quad A, B=0,\cdots, d-1.
\]
Then $(F_{AB})$ is the metric matrix which is positive
Hermitian. Let
\[
F_{AB}=\sum_{C=0}^{d-1}
G_{AC}\bar{G_{BC}}
\]
for a $d\times d$ matrix $(G_{AB})$
and let $(H_{AB})$ be the inverse matrix of $(G_{AB})$. Then it is easy to
see that $\{\sum_{B=0}^{d-1}H_{AB}S_B\}, (A=0,\cdots, d-1)$ forms an
orthonormal basis of $H^0(M,L^m)$. By the definition of 
$S_A (A=0,\cdots, d-1)$, we have
\begin{equation}\label{delta}
\sum_A||\sum_B H_{AB} S_B(x_0)||_{h_m}^2=\sum_A||H_{A0}S_0(x_0)||_{h_m}^2
=\sum_{A=0}^{d-1}|H_{A0}|^2||S_0(x_0)||^2_{h_m}.
\end{equation}
Suppose $(I_{AB})$ is the inverse matrix of $(F_{AB})$. Then 
by the definition of $(H_{AB})$, we have
\begin{equation}\label{eq22}
\sum_{A=0}^{d-1}|H_{A0}|^2=I_{00}.
\end{equation}

By the above discussion, we know that in order to prove Theorem
~\ref{fud}, we just need to estimate the quantity $I_{00}$ and 
$||S_0||_{h_m}$ at $x_0$ for a suitable choice of the basis
$S_0,\cdots,S_{d-1}$.
We use Tian's peak section method in~\cite{T5}
to get the estimates.

We construct peak sections of $L^m$ for $m$ large. Choose a local normal
coordinate $(z_1,\cdots,z_n)$ centered at $x_0$  such that
the Hermitian matrix $(g_{\alpha\bar\beta})$ satisfies
\begin{align}
\begin{split}
& g_{\alpha\bar\beta}(x_0)=\delta_{\alpha\beta}\\
&\frac{\pa^{p_1+\cdots+p_n}g_{\alpha\bar\beta}}{\pa z_1^{p_1}\cdots\pa
z_n^{p_n}}
(x_0)=0
\end{split}
\end{align}
for $\alpha, \beta=1,\cdots,n$ and any nonnegative integers $p_1,\cdots,
p_n$ with $p_1+\cdots+p_n\neq 0$. Such a local
coordinate
system exists and is unique up to an affine transformation. 
This coordinate system is known as the $K$-coordinate. 
See
~\cite{Bo} or~\cite{ru} for
details.

Next we choose a local holomorphic frame $e_L$ of $L$ at $x_0$
such that the local representation function $a$ of the Hermitian metric
$h$ has the properties
\begin{equation}
a(x_0)=1, \frac{\pa^{p_1+\cdots+p_n}a}{\pa z_1^{p_1}\cdots\pa
z_n^{p_n}}(x_0)=0
\end{equation}
for any nonnegative integers $(p_1,\cdots,p_n)$ with $p_1+\cdots+p_n\neq
0$.

Suppose that the local coordinate $(z_1,\cdots,z_n)$ is defined on an
open neighborhood $U$ of $x_0$ in $M$. Define 
the function $|z|$ 
by 
$|z|=\sqrt{|z_1|^2+\cdots+|z_n|^2}$  for $z\in U$.

Let $\Z_+^n$ be the set of $n$-tuple of integers $(p_1,\cdots,p_n)$
such that $p_i\geq 0 (i=1,\cdots, n)$. Let $P=(p_1,\cdots,p_n)$.
Define
\begin{equation}\label{zp}
z^P=z_1^{p_1}\cdots z_n^{p_n}
\end{equation}
and 
\[
p=p_1+\cdots+p_n.
\]

The following lemma 
is  proved in ~\cite{T5} 
using 
the standard $\bar\pa$-estimates (see e.g.~\cite{Hor}).

\begin{lemma}\label{tian}
For  $P=(p_1,\cdots,p_n)\in \Z_+^n$, and an 
integer $p'>p=p_1+\cdots+p_n$, there exists an $m_0>0$ such that for
$m>m_0$, there is a holomorphic global section $S_{P,m}$ in $H^0(M,L^m)$,
satisfying
\[
\int_M||S_{P,m}||^2_{h_m} dV_g=1,\qquad
\int_{M\backslash\{\rr\}}
||S_{P,m}||^2_{h_m} dV_g=\oo{2p'},
\]
and $S_{P,m}$ can be decomposed as 
\[
S_{P,m}=\tilde S_{P,m}+u_{P,m},\qquad (\tilde S_{P,m} \,
\text{and}\,u_{P,m}\,
\text{not necessarily continuous})
\]
such that
\[
\tilde S_{P,m}(x)=\left\{
\begin{array}{ll}
\lambda_Pz^Pe_L^m(1+\oo{2p'}) &x\in\{ \rr\}\\
0                        &x\in M\backslash\{\rr\},
\end{array}
\right.
\]
\[
u_{P,m}(x)=\ooo{2p'}\qquad x\in U,
\]
and
\[
\int_M||u_{P,m}||_{h_m}^2 dV_g=\oo{2p'},
\]
where  
$\oo{2p'}$ denotes a quantity dominated by $C/m^{2p'}$ with the constant
$C$ depending only on $p'$ and the geometry of $M$. Moreover
\[
\lambda_P^{-2}
=\int_\rr |z^P|^2a^mdV_g.
\]
\end{lemma}

\qed

Because of the above lemma, in the rest of this paper, we will use
$S_{P,m}^{p'}$ to denote the peak sections defined above. Furthermore, we
always set
\[
S_0=S_{(0,\cdots,0),m}^{p'}
\]
for $p'$ and $m$ large enough.

We use the notation $|||\cdot|||$
to denote the $L^2$ norm of a section of $H^0(M,L^m)$. That is, if
$T\in H^0(M,L^m)$, then
\[
|||T|||=\sqrt{\int_M||T||^2_{h_m} dV_g}.
\]

The following  lemma~\cite[Lemma 3.2]{ru} is a generalization of
the
lemma
of
Tian~\cite[Lemma 2.2]{T5}. 

\begin{lemma}[Ruan]\label{ruan}
Let $S_P=S_{P,m}^{p'}$ be the section constructed
in the above lemma. Let $T$ be another section of $L^m$.
Near $x_0$, $T=fe_L^m$ for a holomorphic function $f$. When we say $T$'s 
Taylor expansion at $x_0$, we mean the Taylor expansion of $f$ at $x_0$
under the coordinate system $(z_1,\cdots,z_n)$.
\begin{enumerate}
\item If $z^P$ is not in $T$'s Taylor expansion at $x_0$, then
\[
(S_P, T)_{h_m}=\oo{}|||T|||.
\]
\item If $T$ contains no terms $z^Q$, such that $q<p+\sigma$ $(1\leq 
\sigma\leq
p'-n-p$, and
 $\sigma$ is an integer. Recall that
$p=p_1+\cdots+p_n$ and $q=q_1+\cdots+q_n$.) 
in the  Taylor expansion at $x_0$, then
\begin{equation}\label{s1}
(S_P, T)_{h_m}=\oo{1+\sigma/2} |||T|||.
\end{equation}
\end{enumerate}
\end{lemma}

{\bf Proof:}
We only prove 2, since 1 is similar and easier.

By Lemma~\ref{tian},  we have the decomposition
\[
S_{P,m}^{p'}=\tilde  S_{P,m}+u_{P,m}=\tilde S_P+u_P
\]
with $|||u_P|||^2=\oo{2p'}$. Thus
\[
(u_P, T)_{h_m}=\oo{p'}|||T|||=\oo{1+\sigma/2}|||T|||
\]
by the Cauchy inequality. For the $\tilde S_P$ part, we have
\begin{align*}
&\int_{\rr}<\tilde S_P, T>_{h_m} dV_g\\
&\qquad\qquad=
\lambda_P\int_\rr <\tilde S_P', T>_{h_m} 
dV_g+\oo{2p'}\lambda_P|||T|||,
\end{align*}
where $\tilde S_P'=z^Pe_L^m$ for $\rr$ and is zero otherwise.

Let $dV_0$ be the Euclidean volume. i.e.,
\[
dV_0=(\bb)^ndz_1\wedge d\bar z_1\wedge\cdots\wedge
dz_n\wedge d\bar z_n.
\]
We have
\[
dV_g\geq c dV_0,\qquad a^m\geq ce^{-m|z|^2}
\]
on $U$
for a suitable constant $c>0$ and for $m$ large. Thus we have
\[
\int_\rr |z^P|^2 a^mdV_g\geq c^2\int_\rr |z^P|^2 e^{-m|z|^2} dV_0.
\]
By the simple combinatorial identity
\begin{equation}\label{lem1}
\int_{\C^n}|z^P|^2e^{-m|z|^2} dV_0=\frac{P!}{m^{n+p}},
\end{equation}
where $P!=p_1!\cdots p_n!$, we see that 
\[
\int_\rr |z^P|^2 a^mdV_g\geq c_1\frac{1}{m^{n+p}}
\]
for some number $c_1>0$ independent to $m$. Thus by the 
definition of $\lambda_P$,
\begin{equation}\label{x1}
\lambda_P\leq \frac{1}{\sqrt{c_1}} m^{(n+p)/2}.
\end{equation}

In order to prove \eqref{s1},  we then just  need to prove
\[
\lambda_P\int_\rr<\tilde S_P',
T>_{h_m}dV_g=\oo{1+\frac{\sigma}{2}}|||T|||.
\]
Let
\begin{equation}\label{dxi}
\begin{array}{l}
\xi(x)=\log a+|z|^2\\
\eta(x)=\log\det g_{\alpha\bar\beta}.
\end{array}
\end{equation}
If $\rr$, then $|m\xi|\leq\frac{C}{\sqrt{m}}$.
Since
\[
|e^{m\xi}-1- m\xi-\cdots-\frac{1}{(\sigma+4)!}m^{\sigma+4}
\xi^{\sigma+4}|\leq m^{\sigma+5}|\xi|^{\sigma+5}e^{m|\xi|},
\]
we have
\begin{align}\label{p1}
\begin{split}
&\lambda_P|\int_\rr e^{-m\xi+m|z|^2}<\tilde S_P', T>_{h_m}(e^{m\xi}-1-
m\xi\\
&\qquad-\cdots-\frac{1}{(\sigma+4)!}m^{\sigma+4}
\xi^{\sigma+4})e^{-m|z|^2}dV_g|\\
&\leq C\lambda_P\int_\rr|<\tilde S_P', T>_{h_m}|
m^{\sigma+5}|\xi|^{\sigma+5}dV_g\\
&\leq C\lambda_P\sqrt{\int_\rr{|z^P|^{2}m^{2\sigma+10}|\xi|^{2\sigma+10}}
e^{-m|z|^2}dV_g}|||T|||,
\end{split}
\end{align}
where $C$ is a constant independent to $m$. 
Using $|m\xi|\leq\frac{C}{\sqrt{m}}$ again,
we have
\[
\lambda_P\sqrt{\int_\rr{|z^P|^{2}m^{2\sigma+10}|\xi|^{2\sigma+10}
e^{-m|z|^2}dV_g}}=\oo{\sigma+5}=\oo{1+\sigma/2}.
\]
Thus
in order to prove the lemma, we  need only to prove that for any
$k\leq\sigma+4$,
\[
\lambda_P\int_\rr e^{-m\xi+m|z|^2}
<\tilde S_P', T>_{h_m}m^k\xi^ke^{-m|z|^2}dV_g=\oo{1+\frac{\sigma}{2}
}|||T|||.
\]

Let $\xi=\xi_1+\xi_2$ be the decomposition of $\xi$ 
such that $\xi_1$ contains those terms of order less than 
or equal to $4\sigma
+12$
and $\xi_2$ contains those terms of order greater than $4\sigma+12$.
Using
the
similar method as above, we can proved that
\[
\lambda_P\int_\rr <\tilde S_P', T>_{h_m}m^k((\xi_1+\xi_2)^k-\xi_1^k)e^{-m\xi}dV_g
=\oo{1+\frac{\sigma}{2}}|||T|||
\]
because $|(\xi_1+\xi_2)^k-\xi_1^k|\leq C/m^{2\sigma+5}$.
Thus we only need to prove that
\[
\lambda_P\int_\rr <\tilde S_P', T>_{h_m}m^k\xi_1^ke^{-m\xi}dV_g
=\oo{1+\frac{\sigma}{2}}|||T|||
\]
for $k\leq\sigma+4$.
Let
\[
e^\eta=\eta_1+\eta_2,
\]
where $\eta_1$ consists of the terms of order less than
or equal to  $4\sigma+12$ and $\eta_2=e^\eta-\eta_1$. 
Then as above
\[
\lambda_P\int_\rr  <\tilde S_P', T>_{h_m}m^k\xi_1^k\eta_2e^{-m\xi}dV_0
=\oo{1+\frac{\sigma}{2}}|||T|||,
\]
where $dV_0$ is the Euclidean volume form.

It remains to prove that
\[
\lambda_P\int_\rr <\tilde S_P', T>_{h_m}m^k\xi_1^k\eta_1e^{-m\xi}dV_0
=\oo{1+\frac{\sigma}{2}}|||T|||.
\]
Note that $\xi_1^k\eta_1$ is a polynomial of $z$ and $\bar z$ whose
number of terms is bounded by a constant only depending on $\sigma$
and $n$.
Let
\begin{equation}\label{exp}
\xi_1^k \eta_1
=\sum\xi_{IJ}z^I\bar z^J.
\end{equation}
If $|I|-|J|<\sigma$, then by the assumption on $T$,
\[
\int_\rr <\tilde S_P', T>_{h_m}m^k\xi_{IJ}z^I\bar z^Je^{-m\xi}dV_0=0.
\]
On the other hand, 
under the $K$-coordinates, in the expansion of $\xi$, there are no
$z^P$ or $z^P\bar z$ terms. Thus in ~\eqref{exp}, we must have
$|J|\geq 2k$. If $|I|-|J|\geq\sigma\geq 1$, then $|J|\geq 1$. So
\[
|I|+|J|-2k\geq\sigma+2.
\]
Thus 
\begin{align*}
&\lambda_P\int_\rr <\tilde S_P', T>_{h_m}m^k\xi_1^k\eta_1
e^{-m\xi}dV_0
\\
&\leq C\lambda_P\sqrt{
\int_\rr |z^P|^2m^{2k}|z|^{2(2k+\sigma+2)}e^{-m|z|^2}
dV_0}|||T|||.
\end{align*}
The lemma follows from \eqref{x1} and the elementary fact that
\begin{equation}\label{lem2}
\int_{\C^n}
|z^P|^2 |z|^{2q} e^{-m|z|^2} dV_0
=\frac{(n+p+q-1)!P!}{(p+n-1)!m^{n+p+q}}.
\end{equation}

\qed

In the above lemma, if $T$ itself  is  a peak section, then we have a more
accurate result. Before going to the result, let's first define the weight
and the index of a polynomial (resp. monomial, series).

\begin{definition}\label{md1}
Let $R$ be a component of the $i$-th order covariant derivative of the
curvature
tensor, or the Ricci tensor, or the scalar curvature
at a fixed point
where $i\geq 0$. Define the weight
$w(R)$ of $R$ to  be the number $(1+\frac i2)$. For example,
\[
w(\crr ijkl)=w(\css ij)=w(\rho)=1
\]
and
\[
w(R_{i\bar jk\bar l,m})=\frac 32.
\]
The above definition of the weight can be naturally extended to any
monomial of components of the curvature tensor and its derivatives. If a
polynomial
(resp. series) whose any term is of weight $i$, then we call the
polynomial (resp. series) is of  weight $i$.
\end{definition}

More generally, we have the following definition of index.

\begin{definition}\label{md}
Let $f$ be a monomial of the form 
\[
m^\mu gz^P\bar z^Q,
\]
where $\mu$ is an integer or half of an integer, $P,Q$ are multiple index,
$z^P$ and $z^Q$ are defined as in ~\eqref{zp},
and $g$ is a monomial of components of the curvature tensor and its
derivatives at a fixed point. 
In the rest of this paper, a polynomial (resp. series)
is always a polynomial generated by the monomials of the above form.
Define
the index of $f$ by
\[
ind(f)=\mu+w(g)-\frac{p+q}{2},
\]
where $p=p_1+\cdots+p_n$ and $q=q_1+\cdots+q_n$.
We say a polynomial (resp. series) is of index $i$, if all of its 
monomials have the same index $i$. In that case, we say that the 
polynomial (resp. series) is homogeneous.
The polynomials (resp. series) of index $0$ is called
regular. 
\end{definition}

It is easy to check that if $f_1,f_2$ are homogeneous polynomials (resp.
series), then
\[  
ind(f_1f_2)=ind(f_1)+ind(f_2).
\]
The regular polynomials (resp. series) form a ring
under the addition and multiplication.

When a polynomial (resp. monomial, series) contains
no $m$ or $z$'s, the index is the same as the weight.
The following two lemmas give the motivation of the above definitions.

\begin{lemma}\label{n33}
The Taylor expansion  of the function
\[
a^m\det (g_{\alpha\bar \beta})e^{m|z|^2}.
\]
at $x_0$
is a regular series.
\end{lemma}

{\bf Proof:} Consider the Taylor expansion of $\xi$ and $\eta$
in~\eqref{dxi} under the $K$-coordinates. It is not hard to see that the
Taylor expansion of $\xi$
is of index $(-1)$ and the Taylor expansion of $\eta$ is regular. Thus the 
Taylor expansion of $m\xi+\eta$ is regular and so is
the Taylor expansion of
\[
a^m\det g_{\alpha\bar \beta}e^{m|z|^2}=e^{m\xi+\eta}.
\]

\qed

\begin{lemma}\label{lem334}
If 
\[
A_1+\cdots+A_t
\]
is a polynomial of index $i$, then there is a polynomial $B$
of index $i-n$ such that for any $p'>0$,
\[
\int_\rr A_1e^{-m|z|^2}dV_0+\cdots+\int_\rr A_t
e^{-m|z|^2}dV_0-B
=\oo{p'}.
\]
\end{lemma}

{\bf Proof:} Suppose that
\[
A_k=m^\mu gz^P\bar z^Q, \quad 1\leq k\leq t.
\]
If $P\neq Q$, then
\[
\int_\rr A_ke^{-m|z|^2}dV_0=0.
\]
If   $P=Q$, then by \eqref{lem2},
\[
\int_\rr A_k e^{-m|z|^2}dV_0=Cm^{\mu-p-n}g+\oo{p'}
\]
for any $p'>0$
where $C$ is  a constant. Since $A_k$ is of index $i$, we have
\[
\mu+w(g)-p=i.
\]
Thus
\[
ind(Cm^{\mu-p-n}g)=\mu-p-n+w(g)
=i-n.
\]
The lemma is proved.

\qed

\begin{prop}\label{spq}
We have the following expansion for any
$p'>t+2(n+p+q)$,
\[
(S_{P,m}^{p'},S_{Q,m}^{p'})
=\frac{1}{m^\delta}(a_0+\frac{a_1}{m}+\cdots+\frac{a_{t-1}}{m^{t-1}}+\oo{t}),
\]
where $\delta=1$ or $1/2$ and
where all the $a_i$'s are polynomials of the curvature and its derivatives
such that
\[
ind (a_i)=i+\delta.
\]
In particular, the series is regular.
Moreover, all $a_i, (1\leq i\leq t-1)$
can be found by finite steps of algebraic operations from the curvature
and its derivatives.
\end{prop}

{\bf Proof:} The expansion of the function
$z^P\bar z^Q e^{m\xi+\eta}$ has index $-\frac{p+q}{2}$
by Lemma~\ref{n33}. Thus by 
Lemma~\ref{lem334}, there is a polynomial $B_{PQ}$ of 
index $(-\frac{p+q}{2}-n)$ of the form
\[
B_{PQ}=s_0+\frac{s_1}{m}
+\cdots+\frac{s_{t-1}}{m^{t-1}}
\]
 such that
\[
\int_\rr z^P\bar z^Q e^{m\xi+\eta}e^{-m|z|^2}dV_0=B_{PQ}+\oo{p'}.
\]
In particular, $m^{n+p}B_{PP}$ will be regular. Furthermore,
\[
B_{PP}\sim\frac{C_P}{m^{n+p}}
\]
for constant $C_P\neq 0$ by the same argument as in Equation~\eqref{x1}. 
Thus the expansion of 
\[
m^{\frac{n+p}{2}}\sqrt{B_{PP}}
\]
is regular and
$\frac{1}{m^{\frac{n+p}{2}}\sqrt{B_{PP}}}$ expands as a regular series. 
The lemma follows from Lemma~\ref{ruan} and  the fact that
\[
(S_{P,m}^{p'},S_{Q,m}^{p'})=\frac{B_{PQ}}{\sqrt{B_{PP}}\sqrt{B_{QQ}}}
+\oo{n+p+q+t}.
\]

\qed

\section{The Iteration Process}
The main result of this section is Theorem~\ref{thm41}. In order to obtain
the result, we basically use the iteration process in the numerical
analysis for finding the inverse matrix of a given tri-diagonal 
matrix.

\begin{theorem}\label{thm41}
Let $x_0\in M$. 
Suppose $\{S_0, S_1,\cdots, S_{d-1}\}$ is a  basis of
\linebreak $H^0(M,L^m)$
with $S_0(x_0)\neq 0$ and $S_A(x_0)=0$ for $A=1,\cdots,d-1$. Let
$(F_{AB})=((S_A, S_B)_{h_m}) (A,B
=0,\cdots, d-1)$ be the metric matrix. Let
$(I_{AB})$ be the inverse matrix of $(F_{AB})$. Then for any positive
integer $p>0$, we have the expansion
\[  
I_{00}=1+\frac{\sigma_3}{m^3}
+\frac{\sigma_{7/2}}{m^{7/2}}
+\cdots+
\frac{\sigma_{p-1}}{m^{p-1}}+
\frac{\sigma_{((2p-1)/2)}}{m^{(2p-1)/2}}
+\oo{p}.
\]
Furthermore, $\sigma_k (k=3,7/2, 4,\cdots, (2p-1)/2)$
are  polynomials of weight $k$ 
of  
the curvature and the derivatives of the curvature at $x_0$.
\end{theorem}

\begin{rem}
Although not needed, a more careful analysis will show that
$\sigma_{(k/2)}=0$ for all odd $k$'s. 
\end{rem}

Before proving the theorem, we need some algebraic preparation.

\begin{definition}
We say $M=\{M(m)\}$ is a sequence of $s\times s$ block matrices with
block
number $t\in \Z$, if for each $m$,
\[
M=M(m)=
\begin{pmatrix}
M_{11}(m)  & \cdots & M_{1t}(m)\\
\vdots  &\ddots  & \vdots\\
M_{t1}(m) &\cdots & M_{tt}(m)
\end{pmatrix}
\]
such that for $1\leq i,j\leq t$, $M_{ij}$ is a
$\sigma(i)\times\sigma(j)$
matrix and
\[
\sum_{i=1}^t\sigma(i)=s,
\]
where $\sigma:\{1,\cdots, t\}\rightarrow \Z_+$ assigns each number
in $\{1,\cdots,t\}$ a positive integer. We say that
$\{M(m)\}$ is of type $A(p)$ for a positive integer $p$, if for any entry
$s$ of the matrix $M$ ,we have
\begin{enumerate}
\item
If $s$ is a diagonal entry of $M_{ii} (1\leq i\leq t)$ ,
then we have the the following
Taylor expansion
\[
s=1+\frac{s_1}{m}+\cdots+\frac{s_{p-1}}{m^{p-1}}+\oo{p}.
\]
\item
If $s$ is not a diagonal entry of $M_{ii} (1\leq i\leq t)$, then we have
the Taylor expansion
\[
s=\frac{1}{m^\delta}(s_0+\frac{s_1}{m}+\cdots+\frac{s_{p-1}}{m^{p-1}}
+\oo{p}),
\]
where $\delta$ is equal to 1 or $\frac 32$;
\item
If $s$ is an entry of the matrix $M_{ij}$ for which $|i-j|=1$, then
$s=\oo{\frac 32}$. In addition, if $i\neq t$ or $j\neq t$, then
we have the Taylor expansion
\[
s=\frac{1}{m^\delta}(s_0+\frac{s_1}{m}+\cdots+\frac{s_{p-1}}{m^{p-1}}
+\oo{p}),
\]
where $\delta$ is equal to 1 or $\frac 32$;
\item
If $s$ is an entry of $M_{ij}$ for which $|i-j|>1$, then
\[
s=\oo{p}.
\]
\end{enumerate}
The set of all quantities $(s_1/m,\cdots,s_{p-1}/m^{p-1})$,
or $(\frac{s_0}{m^\delta},
\frac{s_1}{m^{1+\delta}},\cdots, \frac{s_{p-1}}{m^{p+\delta-1}})$ for
 $s$ running from all the entries 
of $M_{ij}$ where $|i-j|\leq 1$ and $i\neq t$ or $j\neq t$ 
are called
the Taylor Data of order $p$.
\end{definition}

\begin{rem}
Since $M_{ij}=\oo{p}$ for $|i-j|>1$, it can be treated as zero when we are
only interested in the expansion of order up to $p-1$. A matrix whose
entries $M_{ij}=0$ for $|i-j|>1$ is called a tri-diagonal matrix. For such
a matrix, we have a simple iteration process for finding its inverse matrix.
\end{rem}

The following proposition is a modification of the iteration
process in the numerical analysis for finding the inverse matrix
of a given tri-diagonal matrix.

\begin{prop}\label{prop41} 
Let $M=M(m)$ be a sequence of $s\times s$ block matrices with block number
$t=p+1$ and be of type $A(p)$.
We further assume that $M=M(m)$ is Hermitian positive. 
Let
\[
N=N(m)=
\begin{pmatrix}
N_{11}  & \cdots  & N_{1t}\\
\vdots  & \ddots  &\vdots\\
N_{t1} &\cdots  &N_{tt}
\end{pmatrix}
\]
be the inverse block matrix of $M(m)$. 
By the inverse block matrix, we mean the inverse matrix of the original matrix with the 
same block partition as that of the original block matrix.
We have the following asymptotic expansion:
\begin{equation}\label{p41}
N_{11}=N_{11}^{(0)}+
\frac{N_{11}^{(1/2)}}{m^{1/2}}+
\frac{N_{11}^{(1)}}{m}
+\cdots+
\frac{N_{11}^{(p-1)}}{m^{p-1}}
+\frac{N_{11}^{((2p-1)/2)}}{m^{(2p-1)/2}}+\oo{p}.
\end{equation}
Furthermore, all entries of $N_{11}^{\alpha}/m^\alpha,(0\leq\alpha\leq
(2p-1)/2)$ 
are polynomials of Taylor Data of order $p$ of $M=M(m)$.
\end{prop}

We need the following elementary
lemma.

\begin{lemma}\label{ele}
Suppose
\[
M=
\begin{pmatrix}
M_{11}  &   M_{12}\\
M_{21}  &   M_{22}
\end{pmatrix}
\]
is an invertible block  matrix which is Hermitian positive.
$\bar{M_{11}^T}=M_{11}$,
$\bar{M_{22}^T}=M_{22}$ and $\bar{M_{12}^T}=M_{21}$. Let
\[
N=
\begin{pmatrix}
N_{11}  &   N_{12}\\
N_{21}  &   N_{22}
\end{pmatrix}
\]
be the inverse  block matrix of $M$. Then
\[
N_{11}=M_{11}^{-1}
+M_{11}^{-1}M_{12}(M_{22}-M_{21}M_{11}^{-1}M_{12})^{-1}
M_{21}M_{11}^{-1}.
\]
\end{lemma}

The proof of the lemma is elementary and is omitted.

\qed

{\bf Proof of the Proposition:} Suppose $p=1$. Then by 
Lemma~\ref{ele}, we have
\[
N_{11}=M_{11}^{-1}+M_{11}^{-1}M_{12}
(M_{22}-M_{21}M_{11}^{-1}M_{12})^{-1}
M_{21}M_{11}^{-1}.
\]
Since $M_{12}=\oo{}$, we see that
\[
N_{11}=M_{11}^{-1}+\oo{}=E(\sigma(1))+\oo{},
\]
where $E(\sigma(1))$ is the $\sigma(1)\times\sigma(1)$
unit matrix.
Thus the proposition is true in the case $p=1$.

Assuming  that when $p=k$, the proposition is true. Let $p=k+1$. 
Using Lemma \ref{ele}, we have
\[
N_{11}=M_{11}^{-1}+M_{11}^{-1}(M_{12}\,M_{13}\,\cdots\,M_{1(k+2)})
\tilde M^{-1}(M_{21}\,\cdots\,M_{(k+2)1})^TM_{11}^{-1},
\]
where by Lemma~\ref{ele},
\begin{align*}
&\tilde M
=
\begin{pmatrix}
M_{22} & \cdots & M_{2(k+2)}\\
\vdots & \ddots & \vdots\\
M_{(k+2)2} &\cdots& M_{(k+2)(k+2)}
\end{pmatrix}\\
&\qquad\qquad-
\begin{pmatrix}
M_{21}\\
\vdots\\
M_{(k+2)1}
\end{pmatrix}
M_{11}^{-1}
\begin{pmatrix}
M_{12}& \cdots &M_{1(k+2)}
\end{pmatrix}.
\end{align*}
By the assumption, $M=M(m)$ is a sequence of matrices with the block
number $k+2$. Since $M_{r1}=\oo{k+1} (r>2)$, it is easy to see that
$\tilde M=\tilde M(m)$ is also a sequence of
block matrix with the block number
$k+1$ and is of type $A(k)$.
Furthermore, we have
\begin{equation}\label{td}
N_{11}=M_{11}^{-1}+M_{11}^{-1}M_{12}\tilde N_{22}'M_{21}M_{11}^{-1}
+\oo{k+1},
\end{equation}
where
\[
\tilde N=
\begin{pmatrix}
\tilde N_{22}' &\cdots& \tilde N_{2(k+2)}'\\
\vdots &\ddots&\vdots\\
\tilde N_{(k+2)2}'&\cdots&\tilde N_{(k+2)(k+2)}'
\end{pmatrix}
\]
is the inverse matrix of $\tilde M$. Then by the induction
assumption, we have
\[
\tilde N_{22}'
=\tilde N^{'(0)}_{22}
+\frac{\tilde N_{22}^{'(1/2)}}{m^{1/2}}
+\frac{\tilde N_{22}^{'(1)}}{m}
+\cdots+
\frac{\tilde N_{22}^{'(k-1)}}{m^{k-1}}+
\frac{\tilde N_{22}^{'((2k-1)/2)}}{m^{(2k-1)/2}}+
\oo{k},
\]
where $~\tilde N_{22}^{'(i/2)}/m^{i/2} (1\leq i\leq (2k-1)/2)$ are
polynomials of
Taylor Data
of order $k$ of $\tilde M=\tilde M(m)$. Thus 
$~\tilde N_{22}^{'(i/2)}/m^{i/2} (1\leq i\leq (2k-1)/2)$
must be polynomials of Taylor Data of order $(k+1)$
of $M=M(m)$. On the other
hand, all the entries of the matrix $M_{11}^{-1}$ are polynomials of 
Taylor Data
of order $(k+1)$ of $M=M(m)$. Since $M_{12}=\oo{}$, the entries of
$N_{11}$ are
polynomials of Taylor Data
of order $(k+1)$ of $M=M(m)$ by ~\eqref{td}.
The proposition is proved.

{\bf Proof of  Theorem~\ref{thm41}:}
For a multiple indices, define $|P|=p_1+\cdots+p_n$.
Suppose 
\[
V_k=\{S\in H^0(M,L^m)|D^QS(x_0)=0\,\text{for}\, |Q|\leq k\}
\]
for $k=1,2,\cdots$, where $Q\in \Z_+^n$ is a multiple indices, and
$D$ is a covariant derivative on the bundle $L^m$. 
$V_k=\{0\}$ for $k$ sufficiently large.
For fixed $p$, let $p'=n+8p(p-1)$. Suppose that $m$ is large enough such
that $H^0(M,L^m)$ is spanned by the $S_{P,m}^{p'}$'s for the
multiple indices $|P|\leq 2p(p-1)$ and $V_{2p(p-1)}$. Let $r=d-\dim
V_{2p(p-1)}$.
Then $r$ only depends on $p$ and $n$.
Let 
$T_1,\cdots,T_{d-r}$ be an orthnormal basis
of $V_{2p(p-1)}$
such that 
\[
(S_{P,m}^{p'},T_\alpha)=0
\]
for $|P|\leq 2p(p-1)$ and $\alpha>r$. Let $s(k)=\dim V_k$
for $k\in \Z$. 
For any $1\leq i,j \leq p$, let
$M_{ij}$ be the matrix formed by $(S_{P,m}^{p'},S_{Q,m}^{p'})$ where
$2p(i-2)\leq|P|\leq 2p(i-1)$ and $2p(j-2)\leq|Q|\leq 2p(j-1)$.
Furthermore,
define $M_{i(p+1)}$ to be the matrix whose entries are $(S_P,T_\alpha)$
for $2p(i-2)\leq|P|\leq 2p(i-1)$ and $1\leq\alpha\leq r$. Define
$M_{(p+1)i}$ to be the complex conjugate of  $M_{i(p+1)}$. Finally, define
$M_{(p+1)(p+1)}$ to be the $r\times r$ unit matrix $E(r)$. Then it is easy
to check that $M=(M_{ij})$ is 
a sequence of block matrices of type $A(p)$ with
the block number $p+1$ by using Lemma ~\ref{ruan}
and Proposition~\ref{spq}.

Define an order $\geq$ on the multiple indices $P$ 
as follows:  $P\geq Q$, if
\begin{enumerate}
\item $|P|>|Q|$ or;
\item $|P|=|Q|$ and $p_j=q_j$ but $p_{j+1}>q_{j+1}$ for some $0\leq j
\leq n$.
\end{enumerate}
Using this order, there is a one-one order preserving correspondence $\kappa$
between $\{0,\cdots, r-1\}$ and $\{P||P|\leq 2p(p-1)\}$.

Define
\[
S_A=\left\{
\begin{array}{ll}
S^{p'}_{\kappa(A),m} & A\leq r-1\\
T_{A-r+1}  & A\geq r
\end{array}
\right..
\]

Comparing the matrix $M$ to the metric matrix $F_{AB}=((S_A, S_B)),
(A,B=0,\cdots, d-1)$, by the choice of the basis, we see that
\[
(F_{AB})=
\begin{pmatrix}
M & 0\\
0 & E(d-2r)
\end{pmatrix},
\]
where $E(d-2r)$ is the $(d-2r)\times (d-2r)$ identity matrix.
If  $N=(N_{ij})$ is  the inverse matrix of $M$, then
$N_{11}$ is an $1\times 1$ matrix and 
\[
I_{00}=N_{11}.
\]
So Proposition~\ref{prop41} gives the desired asymptotic expansion
\[
I_{00}=\sigma_0+
\frac{\sigma_{1/2}}{m^{1/2}}
+\frac{\sigma_1}{m}
+\cdots+
\frac{\sigma_{p-1}}{m^{p-1}}+
\frac{\sigma_{((2p-1)/2)}}{m^{(2p-1)/2}}+\oo{p}.
\]
Moreover, Proposition~\ref{prop41} states that $\sigma_k/m^k,
(k=1/2,\cdots,(2p-1)/2)$ are polynomials of the Taylor Data of order $p$
of $M$.
By Proposition~\ref{spq}, the Taylor Data for the inner products are
regular. Thus $\sigma_k/m^k, (1/2\leq k\leq (2p-1)/2)$ are regular or in
other word, 
 $\sigma_k$ is a polynomial of the curvature and its
derivatives of weight $k$ for $k=1/2,\cdots,(2p-1)/2$.

It remains to show that
\[
\sigma_0=1,\,\sigma_{1/2}=\sigma_1=\sigma_{3/2}=\sigma_2=\sigma_{5/2}=0.
\]
This can be seen using the following argument. 
 Let $S_0,S_1,\cdots,S_{d-1}$ be a basis of $H^0(M,L^m)$.
We suppose that $S_A(x_0)=0$ for $A=1,\cdots,d-1$. We also assume that
$(S_0,S_A)=0$ for $A>1$. Let $c=(S_0,S_1)$. Then
$I_{00}=(1-|c|^2)^{-1}$. By Lemma~\ref{ruan}, we see that
$c=\oo{3/2}$. Thus $I_{00}=1+\oo{3}$.

\section{Proof of Theorem~\ref{fud}}
In this section, we prove Theorem~\ref{fud} except that
we postpone the computation of
 the $a_3$ term to the next section. The method of computing the $a_3$
term is the same as that of
$a_j$ for $j=0,1,2$. We put it off to the next section due to the
complexity in the
computation. 

By Theorem~\ref{thm41} and Equation~\eqref{delta}, we just need to
estimate 
\begin{equation}\label{l1}
|\lambda_{(0,\cdots,0)}|^{-2}=\int_\rr a^mdV_g
\end{equation}
to the term $\frac{1}{m^{n+2}}$, from which the first three 
coefficients $a_0$,
$a_1$ and $a_2$  can be calculated.

First we define our notations in the following
equations~\eqref{note1}-~\eqref{note6}.

The curvature tensor is defined as

\begin{equation}\label{note1}
\crr ijkl =\frac{\pa^2 g_{i\bar{j}}}{\pa z_k\pa\bar{z}_l}
-\sum_{p=1}^n\sum_{q=1}^n
g^{p\bar q}\frac{\pa g_{i\bar q}}{\pa z_k}
\frac{\pa g_{p\bar j}}{\pa \bar z_l}
\end{equation}
for $i,j,k,l=1,\cdots,n$. The Ricci curvature is
\begin{equation}\label{note2}
R_{i\bar j}=-\sum_{k,l=1}^ng^{k\bar l}\crr{i}{j}{k}{l}
\end{equation}
for $i,j=1,\cdots,n$,
and the scalar curvature is the trace of the Ricci curvature
\begin{equation}\label{note3}
\rho=\sum_{i,j=1}^ng^{i\bar j}R_{i\bar j}.
\end{equation}

The covariant derivative with respect to $\frac{\pa}{\pa z_p}$
of the curvature tensor is defined as 
\begin{equation}\label{note5}
R_{i\bar jk\bar l,p}=\frac{\pa}{\pa z_p}\crr ijkl
-\sum_{s=1}^n\Gamma_{ip}^s\crr sjkl-\sum_{s=1}^n\Gamma_{kp}^s
\crr ijsp,
\end{equation}
where $\Gamma_{ij}^k=\sum_{s=1}^n
\gtt ks\gss isj$ is the Christoffel symbol.
Higher derivatives are defined in the similar way.
Finally, the Laplace
operator is denoted by $\Delta$:
\begin{equation}\label{note6}
\Delta=\sum_{i=1}^n\sum_{j=1}^ng^{i\bar j}
\frac{\pa^2}{\pa z_i\pa \bar z_j}.
\end{equation}

As will be made obvious, in order to estimate $|\lambda_{(0,\cdots,0)}|^{-2}$
up
to the term $\frac{1}{m^{n+2}}$, we must use the Taylor expansion of
$\log
a$ up to degree 6 and the Taylor expansion of $\,\log\det g_{\alpha\bar
\beta}\,$
up to degree 4. Suppose we have the following expansions at $x_0$:

\begin{gather}\label{ddef}
\begin{split}
\log a=-|z|^2+e_4+e_5+e_6+\ooo{7}\\
\log\det(g_{\alpha\bar\beta})=c_2+c_3+c_4+\ooo{5},
\end{split}
\end{gather}
where $e_4$, $e_5$  and $e_6$ are homogeneous polynomials of
$z$ and $\bar z$ of degree 4,5  and 6, respectively and 
$c_2$, $c_3$, and $c_4$ are homogeneous polynomials of degree
2,3 and 4,
respectively.
Then a straightforward computation gives
\begin{equation}\label{ev}
\left\{
\begin{array}{l}
e_4=-\frac 14 \crr ijklz_iz_k\bar{z}_j\bar z_l\\
e_5=-\frac{1}{12}R_{i\bar jk\bar l,p}z_iz_kz_p\bar z_j\bar z_l
-\frac{1}{12}R_{i\bar jk\bar l,\bar q}z_iz_k\bar z_j\bar z_l\bar z_q\\
\tilde e_6=-\frac{1}{36}(R_{i\bar jk\bar l,p\bar q}
+\crr ispq\crr sjkl +\crr kspq\crr slij\\
\qquad\qquad
+\crr iskq\crr sjpl)
z_iz_kz_p\bar z_j\bar z_l\bar z_q
\end{array}
\right.,
\end{equation}
where $\tilde e_6$ is the $(3,3)$ part of $e_6$, i.e.,
\[
\tilde e_6=\frac{1}{36}\frac{\pa^6\log a}{\pa z_i\pa z_k\pa z_p
\pa\bar z_j\pa\bar z_l\pa\bar z_q}z_iz_kz_p\bar z_j
\bar z_l\bar z_q,
\]
and 
\begin{equation}\label{cv}
\left\{
\begin{array}{l}
c_2=-R_{i\bar j}z_i\bar z_j\\
c_3=-\frac 12R_{i\bar j,k}z_iz_k\bar z_j
-\frac 12R_{i\bar j,\bar l}z_i\bar z_j\bar z_l\\
\tilde c_4=-\frac 14(R_{i\bar j,k\bar l}+\crr isklR_{s\bar j})
z_iz_k\bar z_j\bar z_l
\end{array}
\right.,
\end{equation}
where $\tilde c_4$ is the $(2,2)$ part of $c_4$.
Here $\crr ijkl$, etc denote the value $\crr ijkl(x_0)$.

Considering the function $e^{m(\log a+|z|^2)}e^{\log\det g_{\alpha\bar
\beta}}$, by
~\eqref{ddef}, we
have
\begin{align*}
&e^{m(\log
a+|z|^2)}e^{\log\det
g_{i\bar j}}=e^{m(e_4+e_5+e_6+\ooo{7})}e^{c_2+c_3+c_4+\ooo{5}}\\
&=(1+m(e_4+e_5+e_6)
+\frac 12 m^2(e_4^2)
+O(\cdots))\\
&\qquad\qquad(1+c_2+c_3+c_4+\frac12(c_2^2)+\ooo{5})\\
&=1+m(e_4+e_5+e_6)+\frac 12 m^2(e_4^2)\\
&\qquad+c_2+mc_2e_4
+c_4+\frac 12(c_2^2)
+O(\cdots)),
\end{align*}
where $O(\cdots)$ represents the sum of terms,
each of which is 
less than a constant multiple of $m^\mu|z|^{t}$ for some $\mu$ and $t$
such that $t-2\mu>4$. Those terms will not affect the value of $a_i$,
$i=0, 1, 2$ and can be omitted.

Let $\ph$ be a function on a neighborhood of 
the original point of $\C^n$. For large $m$, define
\begin{equation}\label{k1}
K(\ph)=\int_{\rr}\ph e^{-m|z|^2} dV_0.
\end{equation}
Since 
the functional $K$ is an integration against a symmetric domain,
we have $K(e_5)=0$.
Thus
\begin{align}\label{lamb}
\begin{split}
&|\lambda_{(0,\cdots,0)}|^{-2}=K(e^{m(|z|^2+\log a)}e^{\log\det
g_{i\bar j}})\\
&=K(1)+mK(e_4)+mK(e_6)
+\frac 12 m^2K(e_4^2)\\
&+K(c_2)+mK(c_2e_4)+K(c_4)+\frac 12 K(c_2^2)+\oo{5/2}.
\end{split}
\end{align}

We need the following  combinatorial lemma and its corollary 
which greatly simplified our computation
in this paper.
In fact, it makes our computation feasible.

Suppose $t>0$ is an integer. A function $A$ on $
\{1,\cdots,n\}^{p}\times \{1,\cdots,n\}^{p}$ is
called symmetric, if
\[
A(\sigma(I), \eta(J))=A(I,J)
\]
where $I,J\in \{1,\cdots,n\}^p$ and $\sigma, \eta\in\Sigma$, the transformation
group
of
$\{1,\cdots,n\}$.

For the sake of simplicity, we set $A_{I,\bar J}=A(I,J)$.

\begin{lemma}\label{lem33} Let
$A$ be a symmetric function on $
\{1,\cdots,n\}^{p}\times \{1,\cdots,n\}^{p}$.
Then for any $p'>0$,
\begin{align*}
&\sum_{I,J}\int_\rr
A_{I,\bar J}z_{i_1}\cdots z_{i_p}\bar{z_{j_1}\cdots z_{j_p}}
|z|^{2q} e^{-m|z|^2} dV_0\\
&=(\sum_I A_{I,\bar I}) \frac{p! (n+p+q-1)!}{(p+n-1)!m^{n+p+q}}
+\oo{p'},
\end{align*}
where $I=(i_1,\cdots,i_p)$, $J=(j_1,\cdots,j_p)$ and $
1\leq i_1,\cdots,i_p,j_1,\cdots,j_p\leq n$.
\end{lemma}
{\bf Proof:} Fixing $I=(i_1,\cdots, i_p)$, suppose
\[
z_{i_1}\cdots z_{i_p}=z_1^{p_1}\cdots z_n^{p_n}.
\]
Then
\[
\sum_{i=1}^n p_i=p.
\]

It is easy to see that if $\sigma(J)\neq I$ for any $\sigma\in\Sigma$,
then
\[
\int_\rr z_{i_1}\cdots z_{i_p}\bar{z_{j_1}\cdots z_{j_p}}
|z|^{2q} e^{-m|z|^2} dV_0=0.
\]
On the other hand,
if $\sigma(I)=J$ for some $\sigma\in\Sigma$, then
 the number of $J$'s such that $\sigma(J)=I$ is
$\frac{p!}{p_1!\cdots p_n!}$.
Thus by \eqref{lem2}, we have
\begin{align*}
&\sum_{I,J}\int_{\C^n} A_{I,\bar J}
z_{i_1}\cdots z_{i_p}\bar{z_{j_1}\cdots z_{j_p}}|z|^{2q} e^{-m|z|^2}
dV_0\\
&=\sum_I A_{I,\bar I}\frac{p!}{p_1!\cdots p_n!}
\int_{\C^n}|z_1^{p_1}\cdots z_n^{p_n}|^2|z|^{2q} e^{-m|z|^2} dV_0\\
&=\sum_I A_{I,\bar I}\frac{p!}{p_1!\cdots p_n!}
\frac{(n+p+q-1)!p_1!\cdots p_n!}{(p+n-1)!m^{n+q+p}}
\qquad\text{by~\eqref{lem2}}\\
&=(\sum_I A_{I,\bar I})\frac{p!(n+p+q-1)!}{
(p+n-1)!m^{n+p+q}}.
\end{align*}
The lemma is proved by observing that
\[
e^{-m(\frac{\log m}{\sqrt{m}})^2}=e^{-(\log m)^2}=\oo{p'}
\]
for any $p'>0$.

\qed

Our prototype of function $A$ is the curvature tensor $(\crr ijkl)$,
which is symmetric. However, in most cases, we 
encounter functions which are not symmetric. Thus the following
corollary is  useful.

\begin{cor}\label{cor31}
Let $A$ be a function on $
\{1,\cdots,n\}^{p}\times \{1,\cdots,n\}^{p}$ (not necessarily
symmetric). Then for $p'>0$,
\begin{align*} 
&\sum_{I,J}\int_\rr
A_{I,\bar J}z_{i_1}\cdots z_{i_p}\bar{z_{j_1}\cdots z_{j_p}}
|z|^{2q} e^{-m|z|^2} dV_0\\
&=(\frac{1}{p!}\sum_{I}\sum_{\sigma\in\Sigma}A_{I,\bar{\sigma(I)}})
\frac{p!(n+p+q-1)!}{(p+n-1)!m^{n+p+q}}
+\oo{p'}.
\end{align*}
\end{cor}

{\bf Proof:}
The symmetrization of $A$ is
\[
\tilde
A_{I,\bar J}=\frac{1}{(p!)^2}\sum_{\sigma,\eta\in\Sigma}
A_{\sigma(I),\bar{\eta(J)}}.
\]

Using Lemma~\ref{lem33}, we have

\begin{align*} 
&\sum_{I,J}\int_\rr
A_{I,\bar J}z_{i_1}\cdots z_{i_p}\bar{z_{j_1}\cdots z_{j_p}}
|z|^{2q} e^{-m|z|^2} dV_0\\
&=\sum_{I,J}\int_\rr
\tilde A_{I,\bar J}z_{i_1}\cdots z_{i_p}\bar{z_{j_1}\cdots z_{j_p}}
|z|^{2q} e^{-m|z|^2} dV_0\\
&=\sum_I \tilde A_{I,\bar I}
\frac{p!(n+p+q-1)!}{(p+n-1)!m^{n+p+q}}+\oo{p'}.
\end{align*}

The corollary then  follows from
the elementary fact that

\[
\sum_I\tilde A_{I,\bar I}=\sum_I\frac{1}{(p!)^2}
\sum_{\sigma,\eta\in\Sigma}A_{\sigma(I),\bar{\eta(J)}}
=\frac{1}{p!}\sum_{I}\sum_{\sigma\in\Sigma}A_{\sigma(I),\bar{\sigma(I)}}.
\]

\qed

\begin{prop}\label{prop51}
We have
\begin{gather*}
K(1)=\frac{1}{m^n}+\oo{n+3}\\
mK(e_4)=\frac 12 \rho\frac{1}{m^{n+1}}+\oo{n+3}\\
K(c_2)=-\rho\frac{1}{m^{n+1}}+\oo{n+3}\\
mK(e_6)=-\frac 16(-\Delta\rho+2|Ric|^2+|R|^2)
\frac{1}{m^{n+2}}+\oo{n+3}\\
\frac 12m^2K(e_4^2)
=\frac 18(\rho^2+4|Ric|^2+|R|^2)\frac{1}{m^{n+2}}
+\oo{n+3}\\
mK(c_2e_4)
=-\frac 12(\rho^2+2|Ric|^2)\frac{1}{m^{n+2}}+\oo{n+3}\\
\frac 12K(c_2^2)=\frac 12(\rho^2+|Ric|^2)\frac{1}{m^{n+2}}
+\oo{n+3}\\
K(c_4)=-\frac 12(\Delta\rho-|Ric|^2)\frac{1}{m^{n+2}}+\oo{n+3}.
\end{gather*}
\end{prop}

{\bf Proof:}
We consider $mK(e_4)$ first. Since
\[
e_4=-\frac 14\crr ijkl z_iz_k\bar z_j\bar z_l,
\]
by using Lemma~\ref{lem33}, we have
\[
mK(e_4)=m(-\frac 14
\crr ijij\cdot\frac{2!}{m^{n+2}})+\oo{n+3}
=\frac 12\rho\frac{1}{m^{n+1}}+\oo{n+3}.
\]
Similarly, we can prove the
formulas for $K(1), K(c_2), mK(e_6)$ and $K(c_4)$. Next we consider
$\frac 12 m^2K(e_4^2)$.
We have
\[
e_4^2=\frac{1}{16}\crr ijkl\crr pqrs
z_iz_kz_pz_r\bar z_j\bar z_l\bar z_q\bar z_s.
\]
Using Corollary~\ref{cor31}, we have

\begin{align*}
&\frac 12m^2K(e_4^2)
=\frac{1}{2}\cdot\frac 1{16}m^2\cdot\frac 16(\crr iijj\crr kkll\\
&\qquad\qquad +4\crr
iijk\crr kjll
+\crr ijkl\crr jilk)\frac{4!}{m^{n+4}}+\oo{n+3}\\
&\quad=\frac 18(\rho^2+4|Ric|^2+|R|^2)
\frac{1}{m^{n+2}}+\oo{n+3}.
\end{align*}

The remaining terms can be treated in the similar way.

\qed

{\bf Proof of Theorem~\ref{fud}.}
By Equation~\eqref{delta} and ~\eqref{eq22} we know that the left hand
side of Equation~\eqref{fud1} is equal to
\[
I_{00}|\lambda_{(0,\cdots,0)}|^2.
\]
By Equation~\eqref{l1}, we see that
\[
|\lambda_{(0,\cdots,0)}|^{-2}=\int_\rr e^{m\xi+\eta} e^{-m|z|^2} dV_0,
\]
where $m\xi+\eta=a^m\det g_{\alpha\bar\beta}e^{m|z|^2}$ is a regular
series
by Lemma~\ref{n33}. Thus by Lemma~\ref{lem334}, we have an asymptotic
expansion of $|\lambda_{(0,\cdots,0)}|^{-2}$ of index $(-n)$ whose
each term can be represented by the polynomial of the curvature and its
derivatives.  
Combining this fact to Theorem~\ref{thm41}, we know
 all the $a_i$'s must be
polynomials of the curvature and its derivatives of weight $i$. 
These terms can be obtained by finite many steps of algebraic
operations.
Moreover, by
~\eqref{lamb}
and and Proposition~\ref{prop51}, we see
that
\begin{align}\label{llink}
\begin{split}
&|\lambda_{(0,\cdots,0)}|^{-2}=
\frac{1}{m^n}(1-\frac 12\rho\frac 1m\\
&
\,-\frac{1}{m^2}(\frac
13\Delta\rho+\frac{1}{24}(|R|^2-4|Ric|^2-3\rho^2))+
\oo{5/2}),
\end{split}
\end{align}
from which 
$a_0, a_1, a_2$ can be calculated by the above equation as follows:
\[
\left\{
\begin{array}{l}
a_0=1\\
a_1=\frac 12 \rho\\
a_2=\frac 13\Delta\rho+\frac{1}{24}(|R|^2-4|Ric|^2+3\rho^2).
\end{array}
\right.
\]
This completes the proof of  Theorem~\ref{fud} except for the $a_3$ term.

\qed

\section{Computation of the $a_3$ Term}
In this section, we compute the $a_3$ term. The
first step is to compute $\sigma_3$ in the expansion of Theorem~\ref{thm41}.

\begin{theorem}\label{thm32}
With all the notations as Theorem~\ref{thm41}, we have
\[
\sigma_3=\frac 14|D'\rho|^2,
\]
where $\rho$ is the scalar curvature of the metric $g$
and $|D'\rho|^2=\sum|\frac{\pa\rho}{\pa z_i}|^2$ under local
normal coordinate system.
\end{theorem}

{\bf Proof.}
Let 
\[
V=\{S\in H^0(M,L^m)|S(x_0)=0, DS(x_0)=0\}.
\]
Then $H^0(M,L^m)$ is spanned by  $S_0=S_{(0,\cdots,0),m}^{p'},
S_1=S_{(1,\cdots,0),m}^{p'},
\cdots,S_n=S_{(0,\cdots,1),m}^{p'}$
and $V$, where $\{S_{P,m}^{p'}\}$ are defined in
Section 2 as peak sections.
 Let
$T_1,\cdots, T_{d-n-1}$ be an orthonormal basis of $V$
with $d=\dim H^0(M,L^m)$. The metric
matrix  $(F_{AB})$ can be represented by block matrix
\[
\begin{pmatrix}
1     &   M_{12}      &  M_{13}\\
M_{21}   &   M_{22}   &  M_{23}\\
M_{31}   &  M_{32}    & E
\end{pmatrix},
\]
where $M_{12}=((S_0,S_1),\cdots, (S_0,S_n))$, 
$M_{22}=((S_i,S_j) (1\leq i,j\leq n))$,
$M_{13}=((S_0,T_{\alpha}), 1\leq\alpha\leq d-n-1)$,
$M_{31}=\bar{M_{13}^T})$, 
$M_{23}=((S_i, T_\alpha), 1\leq i\leq n, 1\leq\alpha\leq d-n-1)$, 
$E$ is the unit matrix and $F_{AB}=\bar{F_{BA}}$.

A straightforward computation using Lemma~\ref{ele} shows that
\[
I_{00}=1+
\begin{pmatrix}
M_{12} \,M_{13}
\end{pmatrix}\tilde M^{-1}
\begin{pmatrix}
M_{21}\\
M_{31}
\end{pmatrix},
\]
where
\[
~\tilde M
=
\begin{pmatrix}
M_{22}  & M_{23}\\
M_{32}  & E
\end{pmatrix}
-
\begin{pmatrix}
M_{21}\\
M_{31}
\end{pmatrix}
\begin{pmatrix}
M_{12}  & M_{13}
\end{pmatrix}.
\]
By Lemma~\ref{ruan}, $M_{12}=\oo{3/2}$ and
$M_{13}=\oo{2}$. Thus we have
\begin{equation}\label{izz}
I_{00}=1+M_{12}M_{21}+\oo{7/2}
=1+\sum_{i=1}^n |(S_0,S_i)|^2+\oo{7/2}.
\end{equation}
By the definition of $S_i (0\leq i\leq n)$, we have
\begin{equation}\label{izz1}
(S_0,S_i)={\lambda_{(0,\cdots,0)}
\lambda_{(0,\cdots,1,\cdots,0)}}
\int_\rr\bar z_i a^m dV_g+\oo{2}.
\end{equation}
It is easy to see (cf. ~\cite{T5}) that
\begin{gather}\label{izz2}
\begin{split}
|\lambda_{(0,\cdots,0)}|^{-2}=\frac{1}{m^n}(1+\oo{})\\
|\lambda_{(0,\cdots,1\cdots,0)}|^{-2}
=\frac{1}{m^{n+1}}(1+\oo{}).
\end{split}
\end{gather}
On the other hand, since
\[
\log a=-|z|^2+\frac 14\crr ijkl z_i\bar z_j z_k\bar z_l+\ooo{5},
\]
and since
\[
\log\det g_{i\bar j}=-\css ij z_i\bar z_j+\ooo{3},
\]
we have
\begin{align*}
&\int_\rr \bar z_i a^m dV_g
=-\frac{1}{12}m\int_\rr R_{p\bar qr\bar s, t}\delta_{iu}
z_pz_rz_t\bar z_q\bar z_s\bar z_ue^{-m|z|^2}dV_0\\
&\qquad\quad-\frac 12 \int_\rr R_{p\bar q,r}
\delta_{iu}z_pz_r\bar z_q\bar z_ue^{-m|z|^2}dV_0+\oo{2}.
\end{align*}
where $dV_0$ is the Euclidean volume form.

Thus by Lemma~\ref{lem33}, we have
\begin{equation}\label{izz3}
\int_\rr\bar z_i a^m dV_g=-\frac 12\frac{\rho}{m^{n+2}}(1+\oo{}).
\end{equation}

Using ~\eqref{izz},~\eqref{izz1}, ~\eqref{izz2} and
~\eqref{izz3}, we see that
\[
\sigma_3=\frac 14|D'\rho|^2.
\]

\qed

We now estimate
\[
|\lambda_{(0,\cdots,0)}|^{-2}=\int_\rr a^mdV_g
\]
to the term $\frac{1}{m^{n+3}}$, from which $a_3$  can be calculated.

We define our notations in the following
equations in addition to Equation~\eqref{note1}-~\eqref{note6}.

\begin{equation}\label{note44}
\begin{array}{l}
|R|^2=\sum_{i,j,k,l=1}^n|\crr ijkl|^2\\
|Ric|^2=\sum_{i,j=1}^n|\css ij|^2\\
|D'\rho|^2=\sum_{i=1}^n|\frac{\pa\rho}{\pa z_i}|^2\\
|D'Ric|^2=\sum_{i,j,k=1}^n|R_{i\bar j,k}|^2\\
|D'R|^2=\sum_{i,j,k,l,p=1}^n|R_{i\bar jk\bar l,p}|^2\\
div\,div\,(\rho Ric)=\sum_{i,j=1}^n (\rho \css ji)_{\bar ji}\\
div\,div\,(R,Ric)=\sum_{i,j,k,l=1}^n (\crr ijkl\css ji)_{l\bar k}\\
R(Ric,Ric)=\sum_{i,j,k,l=1}^n\crr ijkl\css ji\css lk\\
Ric(R,R)=\sum_{i,j,k,l,p,q=1}^n\css ij\crr jkpq\crr kiqp\\
\sigma_1(R)=\sum_{i,j,k,l,p,q=1}^n\crr ijkl\crr lkpq\crr qpji\\
\sigma_2(R)=\sum_{i,j,k,l,p,q=1}^n\crr ijkl\crr piqk\crr jplq\\
\sigma_3(Ric)=\sum_{i,j,k=1}^n\css ij\css jk\css ki,
\end{array}
\end{equation}
where ``$,p$'' represents the covariant derivative in the direction 
$\frac{\pa}{\pa z_p}$. 

In order to estimate $|\lambda_{(0,\cdots,0)}|^{-2}$
up
to the term $\frac{1}{m^{n+3}}$, we must use the Taylor expansion of $\log
a$ up to degree 8 and the Taylor expansion of $\log\det g_{i\bar j}$
up to degree 6. Suppose we have the following expansions at $x_0$:

\begin{gather}\label{ddef1}
\begin{split}
\log a=-|z|^2+e_4+e_5+e_6+e_7+e_8+\ooo{9}\\
\log\det(g_{i\bar j})=c_2+c_3+c_4+c_5+c_6+\ooo{7},
\end{split}
\end{gather}
where $e_4$, $e_5$, $e_6$, $e_7$ and $e_8$ are homogeneous polynomials of
$z$ and $\bar z$ of degree 4,5,6,7 and 8, respectively and 
$c_2$, $c_3$, $c_4$ $c_5$ and $c_6$ are homogeneous polynomials of degree
2,3,4,5 and 6,
respectively.

Considering the function $e^{m(\log a+|z|^2)}e^{\log\det g_{\alpha\bar
\beta}}$, by
~\eqref{ddef1}, we
have

\begin{align*}
&\qquad e^{m(\log
a+|z|^2)}e^{\log\det
g_{i\bar j}}\\
&=e^{m(e_4+e_5+e_6+e_7+e_8+\ooo{9})}e^{c_2+c_3+c_4+c_5+c_6
+\ooo{7}}\\
&=1+m(e_4+e_5+e_6+e_7+e_8)+\frac 12 m^2(e_4^2+e_5^2+2e_4e_5+2e_4e_6)\\
&+\frac 16m^3e_4^3
+c_2+mc_2(e_4+e_5+e_6)+\frac 12m^2c_2e_4^2+c_3+mc_3e_5+c_4\\
&+mc_4e_4+c_5+c_6+\frac
12mc_2^2e_4+\frac 12(c_2^2+c_3^2+2c_2c_3+2c_2c_4)\\
&+\frac 16c_2^3
+O(\cdots),
\end{align*}
where $O(\cdots)$ represents the sum of terms,
each of which is 
less than a constant multiple of $m^\mu|z|^{t}$ for some $\mu$ and $t$
such that $t-2\mu>6$. Those terms will not affect the value of $a_3$,
 and can be omitted.

Let $K(\ph)$ be the functional defined in Equation~\eqref{k1}.
We have $K(e_5)=K(e_7)=K(c_2)=K(c_5)=K(c_2e_5)=K(c_2c_3)=0$.
Thus
\begin{align}\label{l-2}
\begin{split}
&\qquad |\lambda_{(0,\cdots,0)}|^{-2}=K(e^{m(|z|^2+\log a)}e^{\log\det
g_{i\bar j}})\\
&=K(1)+mK(e_4)+mK(e_6)+mK(e_8)+\frac 12 m^2K(e_4^2)\\
&+\frac 12
m^2K(e_5^2)+m^2K(e_4e_6)+\frac
16m^3K(e_4^3)+K(c_2)+mK(c_2e_4)\\
&+mK(c_2e_6)+\frac 12 m^2K(c_2e_4^2)
+mK(c_3e_5)+K(c_4)+mK(c_4e_4)\\
&+K(c_6)+\frac 12mK(c_2^2e_4)
+\frac 12 K(c_2^2)
+\frac 12 K(c_3^2)+
K(c_2c_4)\\
&+\frac{1}{6}K(c_2^3)+O(\cdots).
\end{split}
\end{align}

Suppose  $A$ is a homogeneous polynomial
on $\C^n$:
\[
A=\sum_{I,J}A_{I,\bar J}z_I\bar z_J,
\]
where $z_I=z_{i_1}\cdots z_{i_p}$ and $z_J=z_{j_1}\cdots
z_{j_p}$ for $I=(i_1,\cdots,i_p)$ and $J=(j_1,\cdots,j_p)$. Define
\[
L(A)=\frac{1}{p!}\sum_I\sum_{\sigma\in\Sigma}A_{I,\bar{\sigma(I)}}.
\]
Sometimes we also use $L(A_{I,\bar I})$ to denote $L(A)$.
For example,\linebreak $L(\crr ijkl z_iz_k\bar z_j\bar z_l)=
L(\sum\crr iijj)=-\rho$, the scalar curvature.  
By Lemma~\ref{lem33} and Equation~\eqref{l-2}, we have
\begin{align}\label{last}
\begin{split}
&|\lambda_{(0,\cdots,0)}|^{-2}
=\frac{1}{m^n}(1+\frac 1m(2L(e_4)+L(c_2))\\
&+\frac{1}{m^2}(6L(\tilde e_6)+12L(e_4^2)+6L(c_2e_4)
+L(c_2^2)+2L(\tilde c_4))\\
&+\frac{1}{m^3}(L(c_2^3)+3L(c_3^2)+6L(c_2\tilde c_4)
+24L(\tilde c_4e_4)
+12L(c_2^2e_4)+60L(c_2e_4^2)\\
&+120L(e_4^3)+120L(e_4\tilde e_6)
+24L(c_2\tilde e_6)+24L(c_3e_5)+60L(e_5^2)\\
&
+6L(\tilde c_6)+24L(\tilde e_8))+\oo{4}).
\end{split}
\end{align}

\begin{prop}\label{tmp1}
With all the notations as above, we have
\begin{align}\label{f3}
\begin{split}
&L(c_2^3)+3L(c_3^2)+6L(c_2\tilde c_4)
+24L(\tilde c_4e_4)
+12L(c_2^2e_4)+60L(c_2e_4^2)\\
&+120L(e_4^3)+120L(e_4\tilde e_6)
+24L(c_2\tilde e_6)+24L(c_3e_5)+60L(e_5^2)\\
& 
+6L(\tilde c_6)+24L(\tilde e_8)\\
&=-\frac 18\Delta\Delta\rho+\frac 14|D'\rho|^2+\frac
14|D'Ric|^2-\frac{1}{24}
|D'R|^2+\frac 16\rho\Delta\rho+\frac 38\css ij\rho_{j\bar i}\\
&-\frac{1}{48}\rho^3-\frac{1}{12}\rho|Ric|^2+\frac{1}{48}\rho|R|^2-\frac{1}{12}
\sigma_1(R)+\frac{1}{24}\sigma_2(R)+\frac{1}{6}\sigma_3(Ric)\\ 
&+\frac 14R(Ric,Ric).
\end{split}
\end{align}
\end{prop}

We  postpone the proof of this proposition to the Appendix of this paper.

\qed

\begin{prop}\label{g3}
We have
\begin{align*}
&div\,div\,(R,Ric)=-\css ij\rho_{j\bar i}-2|D'Ric|^2
+\crr jilk\crr ij{,k}l\\
&\qquad -R(Ric,Ric)-\sigma_3(Ric)\\
&div\,div\,(\rho Ric)=2|D'\rho|^2+\css ij\rho_{j\bar i}+\rho\Delta\rho\\
&\Delta|R|^2=-2\crr jilk\crr ij{,k}l+2|D'R|^2+4\sigma_1(R)-2\sigma_2(R)
+2Ric(R,R)\\
&\Delta|Ric|^2=2|D'Ric|^2+2\css
ij\rho_{j\bar i}+2R(Ric,Ric)+2\sigma_3(Ric)\\
&\Delta\rho^2=2|D'\rho|^2+2\rho\Delta\rho.
\end{align*}
\end{prop}

{\bf Proof:} A straightforward computation.

\qed

From \eqref{llink},  we know
\begin{align*}
&|\lambda_{(0,\cdots,0)}|^{-2}=\frac{1}{m^n}
(1-\frac 12\rho\frac 1m\\
&-(\frac 13\Delta\rho+\frac{1}{24}(|R|^2-4|Ric|^2-3\rho^2))\frac{1}{m^2}
+\oo{5/2}.
\end{align*}
The term $a_3$ can be computed from
 Proposition~\ref{tmp1}, Proposition~\ref{g3},
Theorem~\ref{thm32}, Equation~\eqref{delta}, ~\eqref{eq22},
~\eqref{last} and the above
expression.

\qed

\begin{example}
Let $M=CP^n$, $L={\mathcal O}(1)$ be the hyperplane bundle. For any $m$,
\[
\sqrt{\frac{(m+n)!}{P!}} z^P
\]
for $P\in\Z_+^n$ with $|P|=m$ form an orthonormal basis of $H^0(M,L^m)$.
Using this we see that
\begin{align*}
&\qquad\sum_{|P|=m}||\sqrt{\frac{(m+n)!}{P!}}
z^P||^2_{h_m}=\frac{(m+n)!}{m!}\\
&=m^n(1+\frac 12n(n+1)\frac 1m
+\frac{1}{24}n(n+1)(n-1)(3n+2)\frac{1}{m^2}\\
&+\frac{1}{48}n^2(n+1)^2(n-1)(n-2)
\frac{1}{m^3}+\oo{4}).
\end{align*}
\end{example}

\begin{cor}
Riemann-Roch Theorem can be recovered
from  Theorem~\ref{fud}, at least asymptotically. Integration
against $M$ on both side of ~\eqref{fud1} gives
\begin{align*}
&\dim H^0(M,L^m)=m^n(vol(M)+\frac 12 c_1(M)\frac 1m\\
&+\frac{1}{12}(c_2+c_1^2)\frac{1}{m^2}+
\frac{1}{24}c_1c_2\frac{1}{m^3}+\oo{4}).
\end{align*}
\end{cor}

\qed

\section{Appendix}
In this Appendix we prove Proposition~\ref{tmp1}. Using the notations as 
in the previous sections, it is splitted into the
following 13 claims.

\begin{claim}
\[
L(c_2^3)=-\frac 16(\rho^3+3\rho|Ric|^2+2\sigma_3(Ric)).
\]
\end{claim}

{\bf Proof:}
\begin{align*}
&L(c_2^3)=-\frac 16(\css ii\css kk\css pp+\css ii\css kp\css pk
+\css ik\css ki\css pp\\
&+\css ik\css kp\css pi+\css ip\css ki\css pk+\css ip\css kk\css pi)\\
&=-\frac 16(\rho^3+3\rho|Ric|^2+2\sigma_3(Ric)).
\end{align*}

\qed

\begin{claim}
\[
3L(c_3^2)=|D'\rho|^2+\frac 12|D'Ric|^2.
\]
\end{claim}

{\bf Proof:}
\begin{align*}
&3L(c_3^2)=\frac 32
L(R_{i\bar i,k}R_{p\bar p,\bar k})\\&
=\frac 12(R_{i\bar i,k}R_{p\bar p,\bar k}
+R_{i\bar p,k}R_{p\bar i,\bar k}
+R_{i\bar k,k}R_{p\bar p,\bar i})\\&
=\frac 12(2|D'\rho|^2+|D'Ric|^2).
\end{align*}

\qed

\begin{claim}
\[
12L(c_2^2e_4)=\frac 14(\rho^3-2R(Ric,Ric)+4\sigma_3(Ric)+5\rho|Ric|^2).
\]
\end{claim}

{\bf Proof:}
\begin{align*}
&12L(c_2^2e_4)=-3L(\crr iikk\css pp\css rr)\\
&=-\frac 14(\crr iikk\css pp\css rr+\crr iikk\css pr\css rp
+\crr iikp\css pk\css rr+\crr iikp\css pr\css rk\\
&+\crr iikr\css pk\css rp+\crr iikr\css pp\css rk
+\crr ikkp\css pi\css rr+\crr ikkp\css pr\css ri\\
&+\crr ikkr\css pi\css rp+\crr ikkr\css pp\css ri
+\crr ipkr\css pi\css rk+\crr ipkr\css pk\css ri)\\
&=\frac 14(\rho^3-2R(Ric,Ric)+4\sigma_3(Ric)+5\rho|Ric|^2).
\end{align*}

\qed

\begin{claim}
\begin{align*}
&6L(c_2\tilde c_4)=
\frac 12\rho\Delta\rho+\css ij\rho_{j\bar i}-\frac 12\rho|Ric|^2
+R(Ric,Ric).
\end{align*}
\end{claim}

{\bf Proof:}
\begin{align*}
&6L(c_2\tilde c_4)
=\frac 32L(\css pp(R_{i\bar i,k\bar k}+\crr irkk\css ri))\\
&=\frac{1}{2}(\css pp(R_{i\bar i,k\bar k}+\crr irkk\css ri)
+2\css pk(\crr ii{,k}p+\crr irkp\css ri))\\
&=\frac{1}{2}(\rho\Delta\rho-\rho|Ric|^2
+2\css pk\rho_{k\bar p}+2R(Ric,Ric)).
\end{align*}

\qed

\begin{claim}
\begin{align*}
&120L(e_4^3)=\frac{1}{48}\rho^3+\frac
14\rho|Ric|^2+\frac{1}{16}\rho|R|^2-\frac 12R(Ric,Ric)\\
&\quad +\frac 12Ric(R,R)-\frac 16{}\sigma_1(R)-\frac
{1}{24}\sigma_2(R)+\frac{1}{3}
\sigma_3(Ric).
\end{align*}
\end{claim}

{\bf Proof:} Suppose
\[
A_{k\bar k_1l\bar l_1p\bar p_1r\bar r_1}
=L(\crr k{k_1}l{l_1}\crr p{p_1}r{r_1}).
\]
Noting that $A_{k\bar k_1l\bar l_1p\bar p_1r\bar r_1}$
is {\it{not}} symmetric, we have

\begin{align}\label{5-1}
\begin{split}
&L(\crr iijjA_{k\bar kl\bar lp\bar pr\bar r})
=\frac{1}{15}(\crr iijj A_{k\bar kl\bar lp\bar pr\bar r}
+8\crr iijkA_{k\bar jl\bar lp\bar pr\bar r}\\
&+2\crr ijklA_{j\bar il\bar kp\bar pr\bar r}
+4\crr ijklA_{j\bar ip\bar pl\bar kr\bar r}).
\end{split}
\end{align}

We compute the above expression term by term as follows:
\begin{align}\label{5-2}
\begin{split}
&A_{k\bar kl\bar lp\bar pr\bar r}=\frac 16
(\crr kkll\crr pprr+4\crr krll\crr pprk+\crr krlp\crr pkrl)\\
&=\frac 16(\rho^2+4|Ric|^2+|R|^2).
\end{split}
\end{align}
\begin{align}\label{5-3}
\begin{split}
&\crr iijkA_{k\bar jl\bar lp\bar pr\bar r}=
-\css jkL(\crr kjll\crr pprr)\\
&=-\frac 16\css jk(\crr kjll\crr pprr+\crr kjlp
\crr plrr+\crr kjlr\crr pprl\\&
+\crr kllp\crr pjrr+\crr kllr\crr pjrp
+\crr kplr\crr pjrl)
\\
&=-\frac 16\rho|Ric|^2+\frac 13R(Ric,Ric)-\frac 13\sigma_3(Ric)
-\frac 16Ric(R,R).
\end{split}
\end{align}
\begin{align}\label{5-4}
\begin{split}
&\crr ijklA_{j\bar il\bar kp\bar pr\bar r}=
\crr ijklL(\crr jilk\crr pprr)\\
&=\frac 16\crr ijkl(\crr jilk\crr pprr+\crr jilp\crr pkrr+\crr jilr\crr
pkrp
\\
&+\crr jklp\crr pirr+\crr jklr\crr pirp+\crr jplr\crr pirk)\\
&=-\frac 16\rho|R|^2-\frac 23Ric(R,R)+\frac 16\sigma_2(R).
\end{split}
\end{align}
\begin{align}\label{5-5}
\begin{split}
&\crr ijklA_{j\bar ip\bar pl\bar kr\bar r}
=\crr ijklL(\crr jipp\crr lkrr)\\
&=\frac 16\crr ijkl(\crr jipp\crr lkrr+\crr jipk\crr lprr+\crr jipr\crr
lprk\\
&+\crr jppk\crr lirr+\crr jppr\crr lirk+\crr jkpr\crr lirp)\\
&=\frac 13R(Ric,Ric)-\frac 13Ric(R,R)+\frac 13\sigma_1(R).
\end{split}
\end{align}
Thus by ~\eqref{5-2} through ~\eqref{5-5}, we have
\begin{align*}
&L(\crr iijjA_{k\bar kl\bar lp\bar pr\bar r})
=\frac{1}{15}(-\frac 16\rho^3-2\rho|Ric|^2-\frac 12\rho|R|^2+4R(Ric,Ric)\\
&-4Ric(R,R)+\frac 43\sigma_1(R)+\frac 13\sigma_2(R)-\frac
83\sigma_3(Ric)).
\end{align*}

\qed

\begin{claim}
\begin{align*}
&24L(\tilde c_4e_4)=-\frac 14\rho\Delta\rho+\frac 14\rho|Ric|^2-\css ip
\rho_{p\bar i}\\
&-R(Ric,Ric)+\frac 14\crr ipkr\crr pi{,r}k+\frac 14Ric(R,R).
\end{align*}
\end{claim}

{\bf Proof:} 
\begin{align*}
&24L(\tilde c_4e_4)
=24L((\tilde c_4)_{i\bar ik\bar k}(e_4)_{p\bar pr\bar r})\\
&=4(c_4)_{i\bar ik\bar k}(e_4)_{p\bar pr\bar r}
+16(\tilde c_4)_{i\bar pk\bar k}(e_4)_{p\bar ir\bar r}
+4(\tilde c_4)_{i\bar pk\bar r}(e_4)_{p\bar ir\bar k}\\
&=\frac 14(\crr iikk(\crr pp{,r}r+\crr p{\alpha}rr\css{\alpha}p)
+\crr ipkk(\crr rr{,p}i+\crr r{\alpha}pi\css{\alpha}r)\\
&+\frac 14\crr ipkr(\crr pi{,r}k+\crr p{\alpha}rk\css{\alpha}i))\\
&=-\frac 14\rho\Delta\rho+\frac 14\rho|Ric|^2-\css ip\rho_{p\bar i}\\
&-R(Ric,Ric)+\frac 14\crr ipkr\crr pi{,r}k
+\frac 14Ric(R,R).
\end{align*}

\qed

\begin{claim}
\begin{align*}
&24L(c_2\tilde e_6)=-\frac 16\rho\Delta\rho-\frac 12\css ri\rho_{i\bar r}
+\frac 13\rho|Ric|^2\\
&+\frac 16\rho|R|^2-R(Ric,Ric)+\frac 12Ric(R,R).
\end{align*}
\end{claim}

{\bf Proof:} 
\begin{align*}
&24L(c_2\tilde e_6)=-24(\css rr(\tilde e_6)_{i\bar ik\bar kp\bar p})\\
&=-6(\css rr (\tilde e_6)_{i\bar ik\bar kp\bar p}
+3\css ri(\tilde e_6)_{k\bar kp\bar pi\bar r})\\
&=\frac 16\rho(-\Delta\rho+2|Ric|^2+|R|^2)\\
&+\frac 12(\css ri(-\rho_{i\bar r}+\crr ksir\crr skpp
+\crr psir\crr spkk+\crr kspr\crr skip)).
\end{align*}

\qed

\begin{claim}
\begin{align*}
&120L(e_4\tilde e_6)
=\frac{1}{12}\rho\Delta\rho+\frac 12\css kr\rho_{r\bar k}-\frac 14\crr
ipkr\crr pi{,r}k\\
&-\frac 16\rho|Ric|^2-\frac{1}{12}\rho|R|^2+R(Ric,Ric)
-\frac 34Ric(R,R)+\frac 12\sigma_1(R).
\end{align*}
\end{claim}

{\bf Proof:} 
\begin{align*}
&L(e_4\tilde e_6)=\frac{1}{10}
((e_4)_{i\bar ik\bar k}(\tilde e_6)_{l\bar lp\bar pr\bar r}
+6(e_4)_{r\bar ij\bar j}(\tilde e_6)_{k\bar kp\bar pi\bar r}
+3(e_4)_{i\bar pk\bar r}(\tilde e_6)_{l\bar lp\bar ir\bar k})\\
&=\frac{1}{1440}(-\rho(-\Delta\rho+2|Ric|^2+|R|^2)
\\
&+6(-\css ri)(-\rho_{r\bar k}+\crr ksir\crr skpp+\crr psir\crr spkk
+\crr kspr\crr skip)\\
&+3\crr ipkr(-\crr pi{,r}k+\crr lsrk\crr slpi
+\crr psrk\crr sill+\crr lspk\crr slri))\\
&=\frac{1}{1440}(\rho\Delta\rho-2\rho|Ric|^2-\rho|R|^2
+6\css kr\rho_{r\bar k}+12R(Ric,Ric)\\
&-6Ric(R,R)-3\crr ipkr\crr pi{,r}k
+3\sigma_1(R)-3Ric(R,R)+3\sigma_1(R)).
\end{align*}

\qed

\begin{claim}
\begin{align*}
&60L(c_2e_4^2)=-\frac 18\rho^3-\rho|Ric|^2-\frac 18\rho|R|^2\\
&+R(Ric,Ric)-\sigma_3(Ric)-\frac 12Ric(R,R).
\end{align*}
\end{claim}

{\bf Proof:} We use the definition of $A_{k\bar kl\bar lp\bar pr\bar r}$
in Claim 5.

\begin{align*}
&L(c_2e_4^2)=-\frac{1}{80}(\css iiA_{k\bar kl\bar lp\bar pr\bar r}
+4\css jkA_{k\bar jl\bar lp\bar pr\bar r})\\
&=-\frac{1}{80}(\frac 16\rho^3+\frac 43\rho|Ric|^2+\frac 16\rho|R|^2\\
&-\frac 43R(Ric,Ric)+\frac 23Ric(R,R)+\frac 43\sigma_3(Ric)).
\end{align*}

\qed

\begin{claim}
\[
24L(c_3e_5)=-|D'\rho|^2-|D'Ric|^2.
\]
\end{claim}

{\bf Proof:}
\begin{align*}
&24L(c_3e_5)=(R_{i\bar i,k}R_{p\bar pr\bar r,\bar k}
+R_{i\bar p,k}R_{p\bar ir\bar r,\bar k})\\
&=-|D'\rho|^2-|D'Ric|^2.
\end{align*}

\qed

\begin{claim}
\[
60L(e_5^2)=\frac{1}{12}(3|D'\rho|^2+6|D'Ric|^2+|D'R|^2).
\]
\end{claim}

{\bf Proof:}
\begin{align*}
&60L(e_5^2)=\frac{120}{144}L(R_{i\bar ij\bar j,k}
R_{p\bar pr\bar r,\bar k})\\
&=\frac 1{12}(3R_{i\bar ij\bar j,k}R_{p\bar pr\bar r,\bar k}
+6R_{i\bar ij\bar p,k}R_{p\bar jr\bar r,\bar k}
+R_{i\bar pj\bar r,k}R_{p\bar ir\bar j,\bar k}).
\end{align*}

\qed

\begin{claim}
\begin{align*}
&6L(\tilde c_6)=-\frac 16\Delta\Delta\rho+\frac 23\css
ij\rho_{j\bar i}
-\frac 13\crr jipk\crr ij{,k}p\\
&+\frac 23|D'Ric|^2-\frac 16Ric(R,R)+\frac 23R(Ric,Ric)+\frac
13\sigma_3(Ric).
\end{align*}
\end{claim}

{\bf Proof:}
Suppose we have the Taylor expansion of the function $\log a$ under the
local $K$-coordinates:
\begin{align*}
&\log a=-|z|^2+(e_4)_{i\bar jk\bar l}z_iz_k\bar z_j\bar z_l
+(e_5^{(1)})_{i\bar jk\bar lp}z_iz_kz_p\bar z_j\bar z_l\\
&+(e_5^{(2)})_{i\bar jk\bar l\bar q}z_iz_k\bar z_j\bar z_l
\bar z_q
+(\tilde e_6)_{i\bar jk\bar lp\bar q}
z_iz_kz_p\bar z_j\bar z_l
\bar z_q+\ooo{7}\\
&+(2,4)+(4,2),
\end{align*}
where $(r,s)\in\N$ represents the sum of terms of the form $az_I\bar z_J$
such that
$|I|=r, |J|=s$. Those terms are irrelevant to the computation of
$a_3$ and need not to be written out explicitly.
Thus we have
\begin{align*}
&g_{i\bar j}=\delta_{ij}-4(e_4)_{i\bar jk\bar l}z_k\bar z_l
-6(e_5^{(1)})_{i\bar jk\bar lp}z_kz_p\bar z_l
-6(e_5^{(2)})_{i\bar jk\bar l\bar q}z_k\bar z_l\bar z_q\\
&-9(e_6)_{i\bar jk\bar lp\bar q}
z_kz_p\bar z_l\bar z_q+\ooo{5}+(1,3)+(3,1).
\end{align*}
Thus the inverse matrix $g^{i\bar j}$ satisfies
\begin{align}\label{6-1}
\begin{split}
&g^{i\bar j}=\delta_{ij}+4(e_4)_{j\bar ik\bar l}z_k\bar z_l
+6(e_5^{(1)})_{j\bar ik\bar lp}z_kz_p\bar z_l
+6(e_5^{(2)})_{j\bar ik\bar l\bar q}z_k\bar z_l\bar z_q\\
&+9(\tilde e_6)_{j\bar ik\bar lp\bar q}
z_kz_p\bar z_l\bar z_q
+16(e_4)_{j\bar sp\bar q}(e_4)_{s\bar ik\bar l}z_p\bar z_qz_k\bar z_l\\
&+\ooo{5}+(1,3)+(3,1).
\end{split}
\end{align}
On the other hand
\begin{align*}
&\log\det g_{i\bar j}=(c_2)_{i\bar j}z_i\bar z_j+(c_3^{(1)})_{i\bar jk}
z_iz_k\bar z_j+(c_3^{(2)})_{i\bar j\bar l}z_i\bar z_j
\bar z_l\\
&+(\tilde c_4)_{i\bar jk\bar l}z_iz_k\bar z_j\bar z_l
+(c_5^{(1)})_{i\bar jk\bar lp}z_iz_kz_p\bar z_j\bar z_l\\
&+(c_5^{(2)})_{i\bar jk\bar l\bar q}
z_iz_k\bar z_j\bar z_l\bar z_q
+(\tilde c_6)_{i\bar jk\bar lp\bar q}z_iz_kz_p\bar z_j\bar z_l\bar
z_q+\ooo{7}\\
&+(1,3)+(3,1)+(1,4)+(4,1)+(1,5)+(5,1)+(2,4)+(4,2).
\end{align*}
Consequently
\begin{align}\label{6-2}
\begin{split}
&-R_{i\bar j}=\pa_i\bar\pa_j\log\det (g_{\alpha\bar\beta})
=(c_2)_{i\bar
j}+2(c_3^{(1)})_{i\bar jk}
z_k+2(c_3^{(2)})_{i\bar j\bar l}
\bar z_l\\
&+4(\tilde c_4)_{i\bar jk\bar l}z_k\bar z_l
+6(c_5^{(1)})_{i\bar jk\bar lp}z_kz_p\bar z_l
+6(c_5^{(2)})_{i\bar jk\bar l\bar q}
z_k\bar z_l\bar z_q\\&
+9(\tilde c_6)_{i\bar jk\bar lp\bar q}z_kz_p\bar z_l\bar
z_q+\ooo{5}+(0,2)+(2,0)+(0,3)+(3,0)\\
&+(0,4)+(4,0)+(1,3)+(3,1).
\end{split}
\end{align}
Using ~\eqref{6-1} and~\eqref{6-2}, we have
\begin{align*}
&-\frac 14\frac{\pa^4\rho}{\pa z^k\pa\bar z^l\pa z^p\pa\bar z^q}
z_kz_p\bar z_l\bar z_q
=9(\tilde c_6)_{i\bar ik\bar lp\bar q}z_kz_p\bar z_l\bar z_q
+16(e_4)_{j\bar ik\bar l}(\tilde c_4)_{i\bar jp\bar q}
z_k\bar z_l z_p\bar z_q\\
&+9(c_2)_{i\bar j}(\tilde e_6)_{j\bar ik\bar lp\bar q}z_k\bar z_lz_p\bar
z_q
+ 12(c_3^{(1)})_{i\bar jp}
(e_5^{(2)})_{j\bar ik\bar l\bar q}z_k\bar z_lz_p\bar z_q\\&
+ 12(c_3^{(2)})_{i\bar j\bar q}(e_5^{(1)})_{j\bar ik\bar l
p}z_kz_p\bar z_l\bar z_q
+16(c_2)_{i\bar j}(e_4)_{j\bar sp\bar q}(e_4)_{s\bar ik\bar l}
z_p\bar z_qz_k\bar z_l.
\end{align*}
Thus
\begin{align*}
&-\frac 14\frac{\pa^4\rho}{\pa z^k\pa \bar z^k\pa z^p\pa \bar z^p}
=9(\tilde c_6)_{i\bar ik\bar kp\bar p}-\frac 34\css ji\rho_{i\bar j}
+\frac 12\crr jipk\crr ij{,k}p-|D'Ric|^2\\
&+\frac 14 Ric(R,R)-R(Ric,Ric)-\frac 12\sigma_3(Ric).
\end{align*}
Since
\[
\Delta\Delta\rho=\frac{\pa^2}{\pa z^k\pa\bar z^k}(g^{p\bar q}\frac
{\pa^2\rho}{\pa z^p\pa\bar z^q})
=\frac{\pa^4\rho}{\pa z^k\pa\bar z^k\pa z^p\pa\bar z^p}
+R_{q\bar p}\rho_{p\bar q}.
\]
We have
\begin{align*}
&\frac 14R_{i\bar j}\rho_{j\bar i}-\frac 14\Delta\Delta\rho
=9(\tilde c_6)_{i\bar ik\bar kp\bar p}-\frac 34R_{i\bar j}\rho_{j\bar i}
+\frac 12\crr jipk\crr ij{,k}p\\&-|D'Ric|^2
+\frac 14Ric(R,R)-R(Ric,Ric)-\frac 12\sigma_3(Ric).
\end{align*}
Thus
\begin{align*}
&(\tilde c_6)_{i\bar ik\bar kp\bar p}=-\frac{1}{36}
\Delta\Delta\rho+\frac 19\css ij\rho_{i\bar j}-
\frac{1}{18}\crr jipk\crr ij{,k}p\\
&+\frac 19|D'Ric|^2-\frac{1}{36}Ric(R,R)+\frac 19R(Ric,Ric)
+\frac 1{18}\sigma_3(Ric).
\end{align*}

\qed

\begin{claim}
\begin{align*}
&
24L(\tilde e_8)=\frac{1}{24}\Delta\Delta\rho-\frac{5}{12}|D'Ric|^2
-\frac 18|D'R|^2-\frac{7}{24}
\css ij\rho_{j\bar i}+\frac 13\crr jikpR_{i\bar j,p\bar k}\\
&+\frac 16Ric(R,R)-\frac{5}{12}R(Ric,Ric)
-\frac{5}{12}\sigma_1(R)+\frac{1}{12}\sigma_2(R)-\frac 16
\sigma_3(Ric).
\end{align*}
\end{claim}

{\bf Proof:} A tedious but  straightforward computation shows that

\begin{align*}
&\frac{\pa^6\log\det g_{i\bar j}}{
\pa z^k\pa\bar z^k\pa z^p\pa\bar z^p\pa z^r\pa\bar z^r}
=\frac{\pa^6 g_{i\bar i}}
{\pa z^p\pa\bar z^p\pa z^k\pa\bar z^k\pa z^r\pa\bar z^r}\\
&-6\crr jikp\frac{\pa^4 g_{i\bar j}}{\pa z^p\pa\bar z^l\pa z^r\pa\bar z^r}
-6|D'Ric|^2-3|D'R|^2+6Ric(R,R)\\
&+2\sigma_1(R)+2\sigma_2(R)
-2\sigma_3(Ric).
\end{align*}

Thus using Claim 12,

\begin{align*}
& 24^2e_8=\Delta\Delta\rho-10|D'Ric|^2-7\css ij\rho_{j\bar i}+8\crr jikp
R_{i\bar j,p\bar k}\\
&-3|D'R|^2+4Ric(R,R)-10R(Ric,Ric)\\
&-10\sigma_1(R)
+2\sigma_2(R)-4\sigma_3(Ric).
\end{align*}

\qed 

Proposition~\ref{tmp1} follows from the above claims.

\bibliographystyle{abbrv}
\bibliography{bib}

\end{document}